\documentclass[MSNbibl,number,citesort,seceqn,dvips]{arxbj}
\usepackage{upgreek}

\aid{0}
\volume{19}
\issue{5A}
\pubyear{2013}
\firstpage{1714}
\lastpage{1749}
\doi{10.3150/12-BEJ427} 

\makeatletter
\newtheorem{them}{Theorem}[section]
\newremark{rem}{Remark}[section]
\newproclaim{ass}{Assumption}[section]
\newproclaim{simpro}{Simulation Procedure}[section]
\newproclaim{imppro}{Implementation Procedure}[section]
\newremark{ex}{Example}[section]
\newtheorem{lem}{Lemma}[section]
\newcommand{\R}{\mathbb{R}}
\newcommand{\f}{\frac}
\newcommand{\no}{\nonumber}
\newcommand{\ep}{\varepsilon}
\newcommand{\lam}{\lambda}
\newcommand{\al}{\alpha}
\newcommand{\si}{\sigma}
\makeatother

\begin{document}
\begin{frontmatter}

\title{Long-range dependent time series specification}
\runtitle{LRD time series specification}

\begin{aug}
\author[1]{\fnms{Jiti} \snm{Gao}\corref{}\thanksref{1}\ead[label=e1]{jiti.gao@monash.edu}},
\author[2]{\fnms{Qiying} \snm{Wang}\thanksref{2}\ead[label=e2]{qiying@maths.usyd.edu.au}} \and
\author[3]{\fnms{Jiying} \snm{Yin}\thanksref{3}\ead[label=e3]{jiying.yin@monash.edu}}
\runauthor{J. Gao, Q. Wang and J. Yin} 
\address[1]{Department of Econometrics and Business Statistics, Monash
University, Caulfield East, Victoria 3145, Australia. \printead{e1}}
\address[2]{School of Mathematics and Statistics, University of
Sydney, Sydney, NSW 2006, Australia.\\ \printead{e2}}
\address[3]{Department of Econometrics and Business Statistics, Monash
University, Caulfield East, Victoria 3145, Australia. \printead{e3}}
\end{aug}

\received{\smonth{12} \syear{2011}}

%
\begin{abstract}
In this paper we propose using a
nonparametric model specification test for parametric time series
with long-range dependence (LRD).
To establish asymptotic distributions of the proposed test
statistic, we develop new central limit theorems for certain
weighted quadratic forms of stationary time series with LRD. To
implement our proposed test in practice, we
develop a computer-intensive parametric bootstrap simulation
procedure for finding simulated critical values. As a result, our
finite-sample studies demonstrate that both the proposed theory and the
simulation procedure work well, and that the proposed test has little
size distortion and reasonable power.
\end{abstract}

%
\begin{keyword}
\kwd{central limit theorem}
\kwd{Gaussian process}
\kwd{linear process}
\kwd{long-range dependence}
\kwd{parametric time series regression}
\kwd{specification testing}
\end{keyword}

\end{frontmatter}

\section{Introduction}\label{sec1}
Starting approximately two decades ago, model specification testing of
short-range dependent stationary time series
has become a very active research field in both econometrics and
statistics. In the meantime, estimation
of long-range dependent stationary time series models has also been
quite active. To the best of our knowledge, however, model
specification of stationary time series with long-range dependence
(LRD) has not been discussed in the literature. This is probably
related to the unavailability of certain central limit theorems for weighted
quadratic forms of stationary time series with LRD.

Specification test statistics based on nonparametric and
semi-parametric techniques for both independent and short-range
dependent cases have been proposed and studied extensively over the
last two decades. Recent studies have shown that some data sets may
display LRD (see, e.g., Beran \cite{be94};
 Cheng and Robinson \cite{chro94}; Robinson
\cite{ro94}; Baillie and King \cite{baki96}; Anh and Heyde \cite{anhe99}; Robinson \cite{ro03}; Gao
\cite{ga07}). In addition, existing studies (e.g., Hidalgo \cite{hi97}; Robinson
\cite{ro97}; Cs\"{o}rg\'{o} and Mielniczuk \cite{csmi99}; Mielniczuk and Wu \cite{miwu04}) have
examined nonparametric regression analysis of data with LRD.

Recent studies (Cheng and Robinson \cite{chro94}; Anh \emph{et al.} \cite{anwogati99};
Mikosch and Starica \cite{mist04}; Gao and Hawthorne \cite{gaha06}) have shown that some
real data in economics, environment, and finance applications with both
LRD and nonlinearities may be modeled by $Y_t =m(X_t) + e_t, 1\leq
t\leq n$, where the regressors $X_t$ are either fixed designs or
stationary short-range dependent variables, $m(\cdot)$ is an unknown
and probably nonlinear trend function, and the errors $e_t$ are
strictly stationary long-range dependent. In addition, the errors may
be allowed for a more general dependent structure of the form
$e_t=G(X_t, Z_t)$ as proposed by Cheng and Robinson \cite{chro94} and
discussed by others, including Mielniczuk and Wu \cite{miwu04}, where $G(\cdot
, \cdot)$ is a suitable function and $\{Z_t\}$ is a possibly
long-range dependent process.

The key findings of these previous studies suggest that to
avoid misrepresenting the mean function or the conditional mean
function of
long-range dependent data, we should let the data speak for themselves
in terms of specifying the true form of the mean function or the conditional
mean function. This is particularly important for data with LRD,
because unnecessary nonlinearity or complexity in mean functions may cause
erroneous LRD.

To address these issues, we propose to model data with possible LRD,
nonlinearity, and non-stationary using a general nonparametric trend model.
Thus, one of the objectives of the present work is to specify the trend by
constructing a nonparametric kernel-based test. Consider a nonlinear time
series model of the form
%
\begin{equation}\label{eq:1.1}
Y_t = m(X_t) + e_t,\qquad t=1,2,\ldots, n,
\end{equation}
where $n$ is the number of observations, $\{X_t\}$ is either a sequence of
fixed designs of $X_t=\frac{t}{n}$ or random regressors, $m(\cdot)$
is an unknown function, and $\{e_t\}$ is a stationary long-range
dependent linear process with $E[e_1]=0$ and $0<E[e_1^2]=\sigma
^2<\infty$. In addition, $\{X_s\}$ and $\{e_t\}$ are assumed to be
independent for all $s, t\geq1$ when $\{X_t\}$ is a sequence of random
regressors. We develop a kernel-based test for the hypotheses
%
\begin{equation}\label{eq:1.2}
H_0\dvt m(x)=m_{\theta_0}(x) \quad\mbox{versus}\quad H_1\dvt m(x) =
m_{\theta_1}(x) + c_n \Delta(x)
\end{equation}
for all $x\in R=(-\infty, \infty)$, where $\theta_0$ and $\theta_1$ are
vectors of
unknown parameters, $m_{\theta}(x)$ is a known parametric function
of $x$ indexed by a vector of unknown parameters, $\theta$, $\Delta
(x)$ is some smooth function, and $\{c_n\}$ is a sequence of real
numbers tending
to 0 when $n\rightarrow\infty$. To ensure that the true model is
identifiable, we also require that for all $\theta_1\neq\theta_2$,
there is a positive constant $\delta_1>0$ such that
$E [m_{\theta_1}(X_1) - m_{\theta_2}(X_1) ]^2 \geq\delta_1>0$ in
the random design case and $\inf_{0\leq x\leq1} |m_{\theta_1}(x) -
m_{\theta_2}(x) | \geq\delta_1>0$ in the fixed design case.

Note that under $H_0$, model
(\ref{eq:1.1}) becomes a parametric model of the form
%
\begin{equation}
Y_t = m_{\theta_0}(X_t) + e_t,
\label{eq:1.3}
\end{equation}
which covers many important cases. For example, model (\ref{eq:1.3}) becomes
a simple linear model with LRD as in (1.1) of Robinson and Hidalgo \cite{rohi97}
when $m_{\theta_0}(X_t) = \alpha_0 + \beta_0 X_t$. For a given set of
long-range
dependent data, the acceptance of $H_0$ suggested by a test statistic may
indicate
that the mean function of the LRD data should be specified parametrically.
In the case of the Nile river data as analyzed by Anh \emph{et al.} \cite{anwogati99},
we would consider using a linear mean function of the form
$m (\frac{t}{n} )=\alpha_0 + \beta_0 \cdot\frac{t}{n}$ if a
suitable test suggested the
acceptance
of $H_0\dvt m(x) = \alpha_0 + \beta_0 x$. Similarly, if a proper test
suggested
accepting a second-order polynomial function of the form
$m (\frac{t}{n} ) = \alpha_0 + \beta_0 (\frac{t}{n} ) + \gamma_0
(\frac{t}{n} )^2$
as the true trend of a financial data set $\{Y_t\}$,
then we would need only to difference $\{Y_t\}$ twice to generate a stationary
set of the
data.

Although there is a long and rich literature in the field of model
specification for time series models with stationarity (see, e.g., Gao
\cite{ga07}; Li and Racine
\cite{lira07}), little work has been reported on parametric
specification testing of $m(\cdot)$ for the case where $X_t$, $e_t$,
or both may be strictly stationary and long-range dependent time
series. To the best of our knowledge, Gao, Wang and Yin \cite{gawayi10} may be
among the first to consider a parametric specification of $m(\cdot)$
for the case where $\{X_t\}$ is a sequence of stationary Gaussian time
series regressors with possible LRD and $\{e_t\}$ is a sequence of
i.i.d. random errors.

This paper is organized as follows.
We present our proposed test for the hypothesis (\ref{eq:1.2}) in
Section~\ref{sec2.1}.
To investigate the proposed specification test, we investigate the
limiting theorems for
the leading term $M_n(h)$ of our test in Section
\ref{sec2.2}, where
%
\begin{equation}
M_n(h)=\sum_{s=1}^n \sum_{t=1, \neq s}^n e_s a_n(X_s, X_t) e_t,
\label{eq:1.4}
\end{equation}
with $a_n(X_s, X_t)= K (\frac{X_s - X_t}{h} )$, in which $K(\cdot)$
is
a probability kernel function and $h$ is a bandwidth parameter.
We mention that
for the case where $\{X_t\}$ is a sequence of either fixed or random regressors
but $\{e_t\}$ is a sequence of long-range dependent errors,
the problem of establishing limiting distributions for $M_n(h)$ is difficult.
Because there is an involvement of $h\rightarrow0$
into the inside of $K(\cdot)$, existing central limit theorems for
U-statistics of long-range dependent processes
(see Hsing and Wu \cite{hswu04}) are not applicable.
Thus, the limit theorems presented in Section~\ref{sec2.2} are interesting and
useful in and of themselves.
In Section~\ref{sec3} we discuss some important extensions and applications
of the theory established in Section~\ref{sec2}. We present a parametric
bootstrap simulation
procedure
as well as some resulting properties for the fixed design situation
in Section~\ref{sec4}. In Section~\ref{sec4} we also provide an example to
demonstrate
how to implement
the proposed test and the bootstrap simulation procedure in practice.
In Section~\ref{sec5} we conclude the paper with some remarks on extensions.
Mathematical details are relegated to Appendices~\ref{appa} and~\ref{appb}. Additional
details are available in Appendix C of Gao and Wang \cite{gawa10}.
Throughout the paper, we use let $a_n\sim b_n$ denote $\lim_{n\to
\infty}a_n/b_n=1$.
\section{Asymptotic theory}\label{sec2}
We propose a test statistic for the hypothesis (\ref
{eq:1.2}) in Section~\ref{sec2.1}. To investigate the proposed test statistic,
we develop some new limiting distributions
of
weighted quadratic forms of dependent processes with LRD in
Section~\ref{sec2.2}.
Their proofs, along with other proofs, are relegated to Appendix~\ref{appa}.
%
\subsection{Model specification test}\label{sec2.1}
Let $K$ be a one-dimensional bounded probability density function and
$h$ be a smoothing bandwidth. When $\{Y_t\}$ is a sequence of long-range
dependent random variables, the conventional kernel estimator of
$m(\cdot)$ is defined by
%
\begin{equation}\label{eq:2.1}
\widehat{m}(x) = \frac{(1/nh) \sum_{t=1}^n K ((x -
X_t)/h )Y_t}
{\widehat{f}(x)},
\end{equation}
where $\widehat{f}(x) = \frac{1}{nh} \sum_{t=1}^n K (\frac{x -
X_t}{h} )$
is the density estimate of the marginal density function,
$f(x)$, of $\{X_t\}$ when $\{X_t\}$ is a sequence of stationary random
regressors. When $X_t=\frac{t}{n}$, $\widehat{f}(x) = \frac{1}{nh}
\sum_{t=1}^n
K (\frac{nx - t}
{nh} )$ is a sequence of functions of $x$.
Various asymptotic properties for $\widehat{m}(x)$ have been studied
in the literature (see, e.g., Cheng and Robinson \cite{chro91}; Robinson \cite{ro97};
Anh \emph{et al.} \cite{anwogati99}).
The nonparametric test statistic for the hypothesis
(\ref{eq:1.2}) might be expected to be related to $\widehat{m}(x)$. However,
as demonstrated in the model specification literature
for both the independent and short-range dependent stationary time
series cases
(see, e.g., Li and Wang \cite{liwa98}; Gao \cite{ga07}; Li and Racine \cite{lira07}), some classes of
nonparametric test statistics have little involvement with
nonparametric estimation.
One of the advantages of such test statistics is that both their
large-sample and finite-sample properties are much less sensitive to
the use of individual nonparametric
estimation as well as the resulting estimation biases.
To test the hypothesis (\ref{eq:1.2}),
we therefore propose a kernel-based test statistic of the form
%
\begin{equation}\label{eq:2.4}
\widehat{M}_n(h) = \sum_{t=1}^n \sum_{s=1,\neq t}^n \widehat{e}_s a_n(X_s,
X_t) \widehat{e}_t
\end{equation}
for the case where $\{X_t\}$ is a sequence of stationary random
regressors, with $a_n(X_s, X_t)= K (\frac{X_s - X_t}{h} )$ and
$\widehat{e}_t = Y_t
- m_{\widetilde{\theta}}(X_t)$, in which $\widetilde{\theta}$ is a
consistent estimator of $\theta_0$ under $H_0$. For the case of fixed-design
mean with LRD errors, we also suggest
the same form of (\ref{eq:2.4}) with $X_t=\frac{t}{n}$. As pointed out
earlier, we chose (\ref{eq:2.4}) over those related to
$\widehat{m}(x)$ mainly because our experience shows that such a form
does not involve biases caused by nonparametric estimation and thus works
well both theoretically and practically.

We need the following assumptions on the error process
$\{e_t\}$, the kernel function $K(\cdot)$, the regressor $\{X_t\}$,
and the regression function
$m_{\theta}(x)$ for the main results of this paper.

\begin{ass}\label{ass2.1}
 \emph{(i)} $\{e_t\}$ is a sequence of strictly
stationary error processes defined by $e_t=\sum_{j=-\infty}^{\infty
}\psi_j\eta_{t-j}$, where the innovations $\{\eta_j\}$ is a sequence
of i.i.d. random variables with $E[\eta_1]=0$, $E[\eta_1^2]=1$ and
$E[\eta_1^6]<\infty$, and the covariance
$\gamma(k)=E[e_t e_{k+t}]=\sum_{j=-\infty}^{\infty}\psi_j\psi
_{j+k}$ satisfies that
$\gamma(0)=\sum_{j=-\infty}^{\infty}\psi_j^2<\infty$ and
$\gamma(k) \sim\eta|k|^{-\alpha}$ as $k\to\infty$, where
$\lambda=(\alpha, \eta)$
with $0<\alpha<1$ and $0<\eta<\infty$ is a vector of unknown parameters.

 \emph{(ii)} In addition, let $\psi_j\ge0$ when $\{X_t\}$ is a sequence of
fixed designs.
\end{ass}
\begin{ass}\label{ass2.2}
 When $\{X_t\}$ is a sequence of i.i.d. random
regressors,
$\{X_t\}$ and $\{e_s\}$ are assumed to be independent for all $s\geq1$
and $t\geq1$, and
the density function $f(x)$ of $\{X_t\}$ is bounded and uniformly continuous.
\end{ass}
\begin{ass}\label{ass2.3}
 Let $\lambda=(\alpha, \eta)$ be defined as in
Assumption~\ref{ass2.1}.
There exist some
$\widetilde{\lambda}=(\widetilde{\alpha}, \widetilde{\eta})$ such
that $\Vert\widetilde{\lambda} - \lambda\Vert=\mathrm{O}_P (w_n^{-1} )$,
where $\{w_n\}$ is a sequence of positive numbers satisfying
$\lim_{n\to\infty} \frac{w_n}{\log n}=\infty$, where $\Vert\cdot\Vert$
denotes the Euclidean norm.\vspace*{-2pt}
\end{ass}
\begin{ass}\label{ass2.4}
 \emph{(i)} $K(\cdot)$ is a bounded and symmetric probability
kernel function over the real line $R$.
 \emph{(ii)} Given $\frac{1}{2}<\alpha<1$, there exists some $0<\beta<\alpha
- \frac{1}{2}$ such that $K(x) =
\mathrm{O} (\frac{1}{1+|x|^{1-\beta}} )$.\vspace*{-2pt}
\end{ass}
\begin{ass}[(Random design)]\label{ass2.5}
 \emph{(i)} Under the null hypothesis
${H}_0$,
$\Vert\widetilde{\theta}-\theta_0\Vert =\mathrm{O}_P(n^{-1/2})$.
 \emph{(ii)} There exists some $\ep_0>0$ such that
$\frac{\partial^2 m_{\theta}(x)}{\partial\theta^2}$ is continuous
in both
$x\in R$ and $\theta\in\Theta_0$, where $\Theta_0=\{\theta\dvt
\Vert\theta-\theta_0\Vert\le\ep_0\}$. In addition,
\[
0<E \biggl[ \biggl\Vert\frac{\partial m_{\theta}(X_1)}{\partial
\theta}\bigg|_{\theta=\theta_0}
\biggr\Vert^2 \biggr]<\infty\quad\mbox{and}\quad 0<E \biggl[ \biggl\Vert\frac{\partial^2 m_{\theta
}(X_1)}{\partial
\theta^2}\bigg|_{\theta=\theta_0}
\biggr\Vert^2 \biggr]<\infty.\vspace*{-2pt}
\]
\end{ass}
\begin{ass}[(Fixed design)]\label{ass2.6}
 \emph{(i)} Under the null hypothesis ${H}_0$,
$\Vert\widetilde{\theta}-\theta_0\Vert =\mathrm{O}_P(n^{-1/2})$.
\emph{(ii)}~There exists some $\ep_0>0$ such that
$\frac{\partial m_{\theta}(x)}{\partial\theta}$ and
$\frac{\partial^2 m_{\theta}(x)}{\partial\theta^2}$ both are
bounded and continuous in
$0\le x\le1$ and $\theta\in\Theta_0$, where $\Theta_0=\{\theta\dvt
\Vert\theta-\theta_0\Vert\le\ep_0\}$.\vspace*{-2pt}
\end{ass}

Assumption~\ref{ass2.1}(i) is quite standard and covers some important cases.
For instance, the case where $\{e_t\}$ is a sequence of Gaussian errors is
included.
Assumption~\ref{ass2.1}(ii) is required to establish Theorems~\ref{them2.2} and~\ref{them2.4} below.
The positivity of $\psi_j$ in Assumption~\ref{ass2.1}(ii) may be replaced by less
restrictive conditions, such as that $\psi_j$ are eventually positive.
Also note that
it is possible to involve a slowly varying function into
the form of $\gamma(k)$; however, because this is not essential to
either our theory or practice,
we use the current Assumption~\ref{ass2.1} throughout this paper.
For alternative conditions on $\{e_t\}$ in this kind of study,
we refer to Cheng and Robinson \cite{chro94}, Robinson \cite{ro97}, and Robinson and
Hidalgo~\cite{rohi97}.

Assumption~\ref{ass2.2} may be relaxed to include the case where $\{X_t\}$ is a
sequence of strictly stationary and $\beta$-mixing random variables.
The corresponding results and their proofs are given in Section 5 and
Appendix C of Gao and Wang \cite{gawa10}.

Assumption~\ref{ass2.3} may be justified for certain $w_n$, like
$w_n=n^{2/5}/\log n$; see the discussion for the construction of
$\widehat{\lambda}$ and Theorem~\ref{them4.2}
in Section~\ref{sec4}. Assumption~\ref{ass2.4}(i) is a standard condition on the kernel function.
To establish Theorems~\ref{them2.2} and~\ref{them2.4}, we also need Assumption~\ref{ass2.4}(ii).
This imposes some restrictions on $\alpha$ and $K(\cdot)$,
but these restrictions are easily verifiable. For instance,
Assumption~\ref{ass2.4}(ii) is satisfied when $K(\cdot)$ is the standard
normal
density function or belongs to a class of probability kernel functions
with compact support. Also note that under Assumption~\ref{ass2.4},
$A_{\alpha}<\infty$, where
%
\begin{eqnarray}\label{ad1}
A_{\alpha}&=&\int_0^{\infty} \int_0^{\infty} \int_0^{\infty}
x^{-\alpha}
y^{-\alpha}
[K(z)K(x+y-z) \no
\\[-8pt]
\\[-8pt]
&&\hphantom{\int_0^{\infty} \int_0^{\infty} \int_0^{\infty}
x^{-\alpha}
y^{-\alpha}
[}{} + K(z-x)K(z-y) ] \,\mathrm{d}x \,\mathrm{d}y \,\mathrm{d}z,\no
\end{eqnarray}
as shown in Lemma~\ref{lema.6}.\eject

The conventional $\sqrt{n}$ rate of convergence assumed in
Assumptions~\ref{ass2.5}(i) and~\ref{ass2.6}(i) is achievable even in this kind of long-range
dependent case (see, e.g., Robinson and Hidalgo \cite{rohi97}; Beran and Ghosh
\cite{begh98}). It is obvious that Assumptions~\ref{ass2.5}(i) and~\ref{ass2.6}(i) imply that
$\Vert\widetilde{\theta}-\theta_0\Vert =\mathrm{o}_P(n^{-\alpha/2})$, which is
crucial and is used in the proofs of Theorems~\ref{them2.1} and~\ref{them2.2}.
Assumption~\ref{ass2.5}(ii) imposes certain moment and smoothing conditions
on both the form $m_{\theta_0}(\cdot)$ and the design $\{X_t\}$.
Assumption~\ref{ass2.6} is used only for the fixed-design case in Theorem~\ref{them2.2}.

We now state the main results of this paper. Theorem~\ref{them2.1} considers the case
where $\{X_t\}$ is a sequence of stationary random regressors.
The case where $X_t=\frac{t}{n}$ is discussed in Theorem~\ref{them2.2}. All notation
is as before.
The proofs of these results are given in Appendix~\ref{appa}.
\begin{them}\label{them2.1}
Suppose that Assumptions~\ref{ass2.1}\emph{(i)}, \ref{ass2.2},
\ref{ass2.3}, \ref{ass2.4}\emph{(i)} and~\ref{ass2.5}
hold.
\begin{longlist}[(ii)]
\item[(i)] If $\lim_{n\rightarrow\infty} n^{2(1-\alpha)} h=0$ and
$\lim_{n\rightarrow\infty} nh=\infty$, then, under $H_0$,
%
\begin{equation}\label{eq:2.5}
\widehat{L}_{1n}(h) \equiv
\frac{\widehat{M}_{n}(h)}{\widehat{\sigma}_{1n}(h)}
\rightarrow_D N(0,1) \qquad\mbox{as } n\rightarrow\infty,
\end{equation}
where $\widehat{\sigma}_{1n}^2(h) =2 n^2 h\int K^2(x) \,\mathrm{d}x
(\frac{1}{n}\sum_{t=1}^n \widehat{f}(X_t) ) (\frac1n
\sum_{t=1}^n\widehat{e}_t^2 )^2$.

\item[(ii)] If $\lim_{n\rightarrow\infty} h=0$ and
$\lim_{n\rightarrow\infty} n^{2(1-\alpha)} h=\infty$,
then, under $H_0$,
%
\begin{eqnarray}\label{eq:2.6}
\widehat{L}_{2n}(h) &\equiv&
\frac{\widehat{M}_{n}(h)}{\widehat{\sigma}_{2n}(h)}
\rightarrow_D \chi^2(1) \qquad\mbox{as } n\rightarrow\infty,
\end{eqnarray}
where $\chi^2(1)$ is the chi-squared distribution with 1 degree of
freedom, and
\[
\widehat{\sigma}_{2n}(h) = \frac{2 n^{2-\widetilde{\alpha}} h
\widetilde{\eta}}
{(1-\widetilde{\alpha})(2-\widetilde{\alpha})}
I_{(0<\widetilde{\alpha}<1)}
\frac{1}{n} \sum_{t=1}^n \widehat{f}(X_t).
\]
\end{longlist}
\end{them}
\begin{them}\label{them2.2}
Suppose that Assumptions~\ref{ass2.1}, \ref{ass2.3}, \ref{ass2.4}
and~\ref{ass2.6} hold.
If $\lim_{n\rightarrow\infty} h=0$, $\lim_{n\to\infty} w_n
h^{1/2}/\log n = \infty$ and $nh\to\infty$ hold,
then, under $H_0$,
%
\begin{equation}\label{eq:2.7}
\widehat{L}_{3n}(h) \equiv \frac{\sum_{t=1}^n \sum_{s=1,\neq
t}^n b_n(s,t)
(\widehat{e}_s \widehat{e}_t - \widehat{\gamma}(s-t) )}
{\widehat{\sigma}_{3n}(h)} \rightarrow_D N(0,1)
\end{equation}
as $n\rightarrow\infty$, where $\widehat{\sigma}_{3n}^2(h) = 8
\widetilde{\eta}^2 n
(nh)^{3-2\widetilde{\alpha}} A_{\widetilde{\alpha}}^*
$ with
\begin{eqnarray*}
A_{\widetilde{\alpha}}^* &=&\int_{1/n}^{n} \int_{1/n}^{n}
\int_{1/n}^{n}
x^{-\widetilde{\alpha}} y^{-\widetilde{\alpha}}
[K(z)K(x+y-z) \no\\
&&\hphantom{\int_{1/n}^{n} \int_{1/n}^{n}
\int_{1/n}^{n}
x^{-\widetilde{\alpha}} y^{-\widetilde{\alpha}}
[}{}+ K(z-x)K(z-y) ] \,\mathrm{d}x \,\mathrm{d}y \,\mathrm{d}z,
\end{eqnarray*}
$b_n(s,t)=K (\frac{s - t}{nh} )$ and
\[
\widehat{\gamma}(k) = \cases{
\displaystyle\frac1n \sum_{i=1}^{n-|k|}\widehat{e}_i\widehat{e}_{i+|k|}
&\quad for $|k|\le(nh)^{1/3}$,\cr
\widetilde{\eta} |k|^{-\widetilde\al} &\quad for $
(nh)^{1/3}<|k|\le n-1$.}
\]
\end{them}
\begin{rem}
(i) As expected, the limiting distributions of
the test statistic $\widehat{M}_n(h)$ (under certain normalization)
for the
hypothesis (\ref{eq:1.2})
depend on the value of $\alpha$ and the choice of the bandwidth $h$.
We require $1/2<\alpha<1$, together with $h^{-1}=\mathrm{o}(n)$ and
$h=\mathrm{o}(n^{-2(1-\alpha)})$, in Theorem~\ref{them2.1}(i). In contrast, Theorem~\ref{them2.1}(ii)
allows $0<\alpha<1$, but we have to restrict $h=\mathrm{o}(1)$ and
$h^{-1}=\mathrm{o}(n^{2(1-\alpha)})$.
These facts imply that to make the limiting distribution of
$\widehat{M}_n(h)$
(under certain normalization) normal,
the conditions $1/2<\alpha<1$ and $h=\mathrm{o}(n^{-2(1-\alpha)})$ are
essentially
necessary for the case where $\{X_t\}$ is a sequence of stationary random
regressors.

(ii) For the fixed-design case of $X_t=\frac{t}{n}$,
we also require $1/2<\alpha<1$
for the asymptotical normality of $\widehat{M}_n(h)$ (under certain
normalization)
in Theorem~\ref{them2.2}. Furthermore,
the range of the bandwidth $h$ depends on the accuracy of
$\Vert\widetilde{\lambda} - \lambda\Vert=\mathrm{O}_P(w_n^{-1})$. Because it may be
justified
that
$w_n=n^{2/5}/\log n$ (see Theorem~\ref{them4.2}), the optimal bandwidth
$h\sim C n^{-1/(4+\al)}$ in theory
is included in Theorem~\ref{them2.2}.
We also mention that the $\widehat\gamma(k)$ defined in Theorem~\ref{them2.2}
provides a consistent estimate of $\gamma(k)$ for each fixed $k$, but
it is not possible to replace
$\widehat{\gamma}(k)=\widetilde\eta|k|^{-\widetilde\al}$
by
$\widehat{\gamma}(k) =\frac{1}{n} \sum_{i=1}^{n-|k|} \widehat{e}_i
\widehat{e}_{i+|k|}$ when
$(nh)^{1/3}<k\le n-1$,
because (\ref{k7}) in the proof of Theorem~\ref{them2.2} below is not true in
the latter case.
The limiting distribution of $\widehat{M}_n(h)$ (under certain
normalization)
for the case of $0<\alpha\le1/2$ in the fixed-design case is an open problem.

(iii) Our experience suggests that it should be possible to establish
some corresponding results of Theorem~\ref{them2.1} for the case where both $\{
X_t\}$ and $\{e_t\}$ exhibit LRD. In the general case, existing studies
(Hidalgo \cite{hi97}; Cs\"{o}rg\'{o} and Mielniczuk \cite{csmi99}; Mielniczuk and Wu
\cite{miwu04}; Guo and Koul \cite{guko07,guko08}; Zhao and Wu \cite{zhwu08}; Kulik and Wichelhaus
\cite{kuwi11}) in nonparametric estimation have already shown that although
similar techniques may be used to establish and prove such
corresponding results, the corresponding conditions and proofs are much
more technical than those involved in the present work. Our preliminary
studies show that we may establish some corresponding results as in
Theorem~\ref{them2.1} for the case where $e_t=\sigma(X_t) \varepsilon_t$, in which
$\{X_t\}$ is a sequence of strictly stationary regressors with LRD and
$\{\varepsilon_t\}$ is a sequence of strictly stationary short-range
dependent errors. To present the key ideas, we focus on the main cases
in Theorems~\ref{them2.1} and~\ref{them2.2}. Other cases are left for future research.
\end{rem}

To motivate the need to establish some limit theorems
for general quadratic forms of dependent processes with LRD, we observe that
%
\begin{equation}
\widehat{M}_n(h) = \sum_{t=1}^n \sum_{s=1,\neq t}^n \widehat{e}_s
a_n(X_s, X_t) \widehat{e}_t
=M_n(h) + 2 R_{1n}(h)+R_{2n}(h),
\label{eq:2.8}
\end{equation}
where $M_n(h) = \sum_{t=1}^n \sum_{s=1,\neq t}^n e_s a_n(X_s, X_t)
e_t $,
\begin{eqnarray*}
R_{1n}(h)&=& \sum_{t=1}^n \sum_{s=1,\neq t}^n a_n(X_s, X_t) e_s
\bigl(m(X_t)
- m_{\widetilde{\theta}}(X_t) \bigr),\\
R_{2n}(h) &=&
\sum_{t=1}^n \sum_{s=1,\neq t}^n a_n(X_s, X_t) \bigl(m(X_s) -
m_{\widetilde{\theta}}(X_s) \bigr) \bigl(m(X_t) -
m_{\widetilde{\theta}}(X_t) \bigr).
\end{eqnarray*}

We show in Appendix~\ref{appa} that $2R_{1n}(h)+R_{2n}(h)=\mathrm{o}_P(\sigma_{in}(h))$,
where $\sigma_{1n}(h)$ and $\sigma_{2n}(h)$ are as defined in Theorem
\ref{them2.3} below. Thus, to prove Theorem~\ref{them2.1}, we need to establish limit
theorems for $M_n(h)$, which is a weighted quadratic form of $\{e_t\}$.
Similar arguments also apply for Theorem~\ref{them2.2}. Because existing results
are not directly applicable to proving Theorems~\ref{them2.1} and~\ref{them2.2}, we
establish our own results the the next section.
%
\subsection{Limit theorems for quadratic forms}\label{sec2.2}
Because both the
conditions and
results for random designs differ
from those for fixed designs, we establish the following
two theorems separately. Their proofs are given in Appendix~\ref{appa}.
\begin{them}\label{them2.3}
Suppose that Assumptions~\ref{ass2.1}\emph{(i)} and~\ref{ass2.2}--\ref{ass2.4}
hold.
\begin{longlist}[(ii)]
\item[(i)] If $\lim_{n\rightarrow\infty} n^{2(1-\alpha)} h =0$
and $\lim_{n\rightarrow\infty} nh=\infty$ hold, then
%
\begin{equation}
\label{eq:2.9}
\frac{\sum_{t=1}^n \sum_{s=1,\neq t}^n e_s a_n(X_s, X_t) e_t
}{\sigma_{1n}(h)}\rightarrow_D N(0,1)
\qquad\mbox{as } n\rightarrow\infty,
\end{equation}
where $\sigma_{1n}^2(h) = n^2 h A_{1\alpha}^2$, in which $A_{1\alpha
}^2 =2 \gamma^2(0) \int K^2(x) \,\mathrm{d}x \int f^2(y) \,\mathrm{d}y$, where $f(\cdot)$
denotes the marginal density of $\{X_t\}$.

\item[(ii)] \emph{If $\lim_{n\rightarrow\infty} h=0$ and
$\lim_{n\rightarrow\infty} n^{2(1-\alpha)}
h =\infty$ hold, then
%
\begin{equation}
\label{eq:2.10}
\frac{\sum_{t=1}^n \sum_{s=1,\neq t}^n e_s a_n(X_s, X_t) e_t
}{\sigma_{2n}(h)}\rightarrow_D
\chi^2(1) \qquad\mbox{as } n\rightarrow\infty,
\end{equation}
where $\sigma_{2n}(h) = n^{2-\alpha} h A_{2\alpha}$,
in which $A_{2\alpha} = \frac{2 \eta}{(1-\alpha)(2-\alpha)} (\int
f^2(y) \,\mathrm{d}y )$}.
\end{longlist}
\end{them}
\begin{them}\label{them2.4}
Suppose that Assumptions~\ref{ass2.1}, \ref{ass2.3} and
\ref{ass2.4} hold.
If in addition both $\lim_{n\rightarrow\infty} h=0$ and $\lim
_{n\rightarrow\infty}
nh=\infty$ hold,
then
%
\begin{equation}
\label{eq:2.11}
\frac{\sum_{t=1}^n \sum_{s=1,\neq t}^n b_n(s,t) (e_s e_t -
\gamma(s-t) )}{\sigma_{3n}(h)}
\rightarrow_D N(0,1)\qquad \mbox{as } n\rightarrow\infty,
\end{equation}
where $\sigma_{3n}^2(h) = 8 \eta^2 n (nh)^{3-2\alpha} A_{\alpha}$,
with $A_{\alpha}$ defined as in (\ref{ad1}).
\end{them}

\begin{rem}
(i) Theorem~\ref{them2.3} extends existing limit theorems
for both i.i.d. and short-range dependent cases
(see, e.g., Gao \cite{ga07}, Chapter 3) to the situation
where $\{e_t\}$ is a long-range dependent linear process.
Unlike the previously reported
results, Theorem~\ref{them2.3} shows that when $\{a_n(X_s, X_t)\}$ is a sequence of
random functions,
the limiting distribution of the random
weighted quadratic form can be
either a standard normal distribution or a chi-squared distribution.

(ii) In the situation where $\{a_n(X_s, X_t)\}$ is a sequence of
non-random functions,
related results on quadratic forms of long-range dependent time series
have been given by Davydov \cite{da70}, Fox and Taqqu \cite{fota87,fota85}, Avram \cite{av88},
Giraitis and Surgailis \cite{gisu90}, Giraitis and Taqqu \cite{gita97},
Ho and Hsing \cite{hohs96,hohs97,hohs03}, Hsing and Wu \cite{hswu04}, and others.
Because the weighted coefficients in those results
(see, e.g., Hsing and Wu \cite{hswu04}) are non-random and independent of $n$,
they are not applicable for the establishment of Theorems~\ref{them2.3} and~\ref{them2.4}.
Thus, both Theorems~\ref{them2.3} and~\ref{them2.4} extend and complement various existing results.
\end{rem}
\section{Extensions and applications}\label{sec3}
In this section we show that each of the leading terms of many existing
kernel-based test statistics may be represented by a quadratic form
similar to (\ref{eq:1.4}). Theorems~\ref{them2.1}--\ref{them2.4} show that it is possible
to establish corresponding results based on the long-range dependent
errors for these existing test
statistics. To avoid some repetitious arguments, here we only state
some main steps; further details are available in Section 3 of Gao and
Wang \cite{gawa10}.
%
\subsection{Existing kernel-based tests for conditional
mean}\label{sec3.1}
A very simple idea for constructing a kernel test for $H_0$ is to compare
the $L_2$ distance between a nonparametric kernel estimator of $m(\cdot)$
and a parametric counterpart. Denote the nonparametric estimator of
$m(\cdot)$ by $\widehat{m}(\cdot)$ as in (\ref{eq:2.1}) and the parametric
estimator of $m_{\theta_0}(\cdot)$ by
$\widetilde{m}_{\widetilde{\theta}}(\cdot)$ given by
%
\begin{equation}
\label{eq;2.14}
\widetilde{m}_{\widetilde{\theta}}(x) = \frac{\sum_{t=1}^n
K_h(x-X_t)m_{\widetilde{\theta}}(X_i)}{\sum_{t=1}^n K_h(x-X_t)},
\end{equation}
where $\widetilde{\theta}$ is a consistent estimator of $\theta_0$ as
defined before and $K_h(\cdot) = \frac{1}{h} K (\frac{\cdot}{h} )$.
H\"ardle and Mammen \cite{hama93} proposed a test statistic of the form
%
\begin{equation}
\label{eq:2.15}
L_{n1}(h)  =  n{h} \int\{\widehat{m}(x) -
\widetilde{m}_{\widetilde{\theta}}(x) \}^2 w(x) \,\mathrm{d}x,
\end{equation}
where $w(x)$ is a non-negative weight function.
Recall the model (\ref{eq:1.3}). Under $H_0$, it is readily seen that
%
\begin{eqnarray}\label{r5}
L_{n1}(h)& = & n{h} \int\biggl(\frac{ [\sum_{t=1}^n
K_h(x-X_t) (e_t + m_{\theta_0}(X_t) - m_{\widetilde{\theta}}(X_t) )
]^2}{n^2 \widehat{f}^2(x)} \biggr) w(x) \,\mathrm{d}x\nonumber\\
& = & n {h} \sum_{s=1}^n \sum_{t=1}^n \biggl(\int
\frac{K_h(x-X_s)K_h(x-X_t)}{n^2 \widehat{f}^2(x)} w(x) \,\mathrm{d}x \biggr)
e_s e_t\\
&&{} + \mbox{remainder term},\no
\end{eqnarray}
where the leading term is similar to (\ref{eq:1.4}). In a closely
related study,
Chen and Gao \cite{chga07} constructed a test statistic based on empirical
likelihood ideas. As they showed, the first-order approximation of
their test is asymptotically equivalent to
%
\begin{equation}
\label{eq:2.18}
L_{n2}(h) = nh \int[ \widehat{m}(x)
-\widetilde{m}_{\widetilde{\theta}}(x)
]^2 w(x) \,\mathrm{d}x.
\end{equation}
It can be easily shown that the test statistic $L_{n2}(h)$
has a similar decomposition to (\ref{r5}),
in which the leading term is similar to (\ref{eq:1.4}).
%
\subsection{Testing conditional mean with conditional
variance}\label{sec3.2}
Because the main objective of this paper is to specify the form of
$m(\cdot)$ parametrically, we have assumed that the variance or
conditional variance $\sigma^2$ is an unknown parameter. As can be
seen from (\ref{eq:1.4}), we can replace $e_t$ by $e_t=\sigma(X_t)
\varepsilon_t$, where
$\{\varepsilon_t\}$
is a sequence of long-range dependent linear processes. In this case,
the leading term of $\widehat M_n(h)$ in (\ref{eq:2.4}) becomes
%
\begin{equation}\label{eq:2.20}
L_{n3}(h)=\sum_{s=1}^n \sum_{t=1, \neq s}^n \varepsilon_s \sigma(X_s)
K \biggl(\frac{X_s-X_t}{h} \biggr) \sigma(X_t) \varepsilon_t,
\end{equation}
which is also a quadratic form of $(X_t, \varepsilon_t)$.
%
\subsection{Testing conditional mean in additive form}\label{sec3.3}
When $X_t=(X_{t1}, \ldots, X_{td})$ in (\ref{eq:1.1}) is a vector of
$d$-dimensional designs, we may consider a hypothetical problem of the form
%
\begin{equation}
H_0^{\prime}\dvt m(x) = \sum_{i=1}^d m_{i\theta_0}(x_i) \quad\mbox{versus}\quad H_1^{\prime}\dvt m(x) = \sum_{i=1}^d m_{i\theta_1}(x_i) + c_n \sum_{i=1}^d
\Delta_i(x_i),
\label{eq:1.5}
\end{equation}
where each $m_{i\theta_0}(\cdot)$ is a known function indexed by
$\theta_0$
and $\Delta_i(\cdot)$ is also a known function over~$\R$. Various additive
models have been discussed in the literature (see, e.g., Sperlich,
Tj{\o}stheim and Yang \cite{sptjya02}; Gao, Lu and Tj{\o}stheim \cite{galutj06}; Gao \cite{ga07}).
The construction of $\widehat M_n(h)$ suggests a test statistic of the
form
%
\begin{equation}
\label{eq:6.13}
L_{n4}(h)=\sum_{j=1}^n \sum_{i=1}^n \widetilde{Y}_i \prod_{k=1}^d
K \biggl(\frac{X_{ik} - X_{jk}}{h}\biggr ) \widetilde{Y}_j
\end{equation}
for the hypothesis (\ref{eq:1.5}),
where $\widetilde{Y}_i = (Y_i - \sum_{k=1}^d m_{k \widetilde{\theta}_0}
(X_{ik}) )$. Clearly, $L_{n4}(h)$ also
has a leading term similar to (\ref{eq:1.4}).

As mentioned earlier, some corresponding results of Theorems~\ref{them2.1}--\ref{them2.4}
may be established accordingly for $L_{ni}(h)$, $1\leq i\leq4$.
\section{Simulation procedure and an example of
implementation}\label{sec4}
 In
this section we focus on an implementation procedure for $\widehat
{L}_{3n}(h)$ in the fixed design case. Implementation of $\widehat
{L}_{1n}(h)$ and $\widehat{L}_{2n}(h)$ in the random design case
requires a different asymptotic theory, which involves developing a
different bootstrap simulation procedure. We leave a discussion of this
for future research.
Before we implement $\widehat{L}_{3n}(h)$, we develop a simulation
procedure for the choice of a simulated critical value and then propose
an estimation procedure for $\lambda$ involved in the proposed test.

Explicitly, Section~\ref{sec4.1} establishes a simulation procedure for the
implementation of the test statistic $\widehat L_{3n}(h)$. An
estimation procedure for $\lambda$ is briefly mentioned in Section
\ref{sec4.2}. Section~\ref{sec4.3} presents an example of implementation to check whether
the theory works in practice.
\subsection{Simulation scheme and asymptotic properties}\label{sec4.1} When
constructing the simulation, the covariance structure $\gamma_{\lambda
}(k)= E[e_t e_{t+k}]$ must be replaced by an estimated version, $\gamma
_{\widetilde{\lambda}}(k)$, with $\widetilde{\lambda}=(\widetilde
{\alpha},\widetilde{\eta})$ as consistent estimators of $\lambda
=(\alpha, \eta)$. We assume the existence of $\widetilde{\lambda}$
at the moment, and describe its construction in Section~\ref{sec4.2}.
\begin{simpro}\label{simpro4.1}
 Let $l_r$ ($0<r<1$) be the $1-r$
quantile of the exact finite-sample distribution of $\widehat
{L}_{3n}(h)$. Because $l_{r}$ might not be evaluated in practice, we
suggest replacing it with an approximate $r$-level critical value
$l_{r}^{\ast}$, using the following bootstrap procedure:
\begin{enumerate}
\item
Generate $Y_t^{\ast}$ from $Y_t^{\ast}=m_{\widetilde{\theta}}(X_t)
+ e_t^{\ast}$ for $1\leq t\leq n$, where $X_t=\frac{t}{n}$ and $\{
e_t^{\ast}\}$ is a sequence of stationary LRD errors generated from
$e_t^{\ast}=\sum_{j=-\infty}^{\infty}\psi_j^{\ast} \eta
_{t-j}^{\ast}$, with $\{\psi_j^{\ast}\}$ chosen as $\psi_j^{\ast}
= c(\widetilde{\alpha}) |j|^{-(1+\widetilde{\alpha})/2}$
such that
$\gamma_{\widetilde{\lambda}}(k)=E[e_t^{\ast} e_{k+t}^{\ast
}|\mathcal{W}_n]=\sum_{j=-\infty}^{\infty}\psi_j^{\ast} \psi
_{j+k}^{\ast}\sim\widetilde{\eta} |k|^{-\widetilde{\alpha}}$ as
$k$ large enough and $c(\widetilde{\alpha})$ also chosen such that
$\sum_{j=-\infty}^{\infty} \psi^{\ast2} = \widetilde{\sigma}^2 =
\frac{1}{n}\sum_{t=1}^n (Y_t - m_{\widetilde{\theta}}(X_t) )^2$,
where $\{\eta_s^{\ast}\}$ is a sequence of independent and
identically random variables with $E[\eta_1^{\ast}]=0$, $E[\eta
_1^{\ast2}]=1$, $E[\eta_1^{\ast6}]<\infty$, and $\mathcal{W}_n=(
Y_1, \ldots, Y_n)$.
\item
Use the data set $\{Y_t^{\ast}\dvt 1\leq t\leq n\}$ to obtain an
estimator $\widetilde{\theta}^{\ast}$ of ${\widetilde\theta}$, and
construct a corresponding version, $\widehat{L}_{3n}^*(h)$, of
$\widehat{L}_{3n}(h)$ under $H_0$; that is, $\widehat{L}_{3n}^*(h)$
is computed according to the same formula as $\widehat{L}_{3n}(h)$ but
with $\{Y_t\dvt 1\leq t\leq n\}$ and ${(\theta_0,\widetilde\theta)}$
replaced by $\{Y_t^{\ast}\dvt 1\leq t\leq n\}$ and $(\widetilde\theta
,\widetilde{\theta}^{\ast})$, respectively.
\item
Repeat the foregoing step $M$ times and produce $M$ versions of
$\widehat{L}_{3n}^{\ast}(h)$, denoted by $\widehat{L}_{3n,m}^{\ast
}(h)$ for $m=1,2,\ldots, M$. Use the $M$ values of $\widehat
{L}_{3n,m}^{\ast}(h)$ to construct their empirical distribution
function. Define $P^{\ast}(\cdot)$ as the bootstrap distribution of
$\widehat{L}_{3n}^*(h)$ given\vspace*{1pt} $\mathcal{W}_n$ by $P^{\ast} (\widehat
{L}_{3n}^*(h)\leq x )=P (\widehat{L}_{3n}^*(h)\leq x|\mathcal{W}_n
)$. Let $l_{r}^{\ast}$ ($0<r<1$) satisfy $P^{\ast} (\widehat
{L}_{3n}^*(h)\geq
l_{r}^{\ast} )=r$, and estimate $l_{r}$ by $l_{r}^{\ast}$.\vspace*{-1pt}
\end{enumerate}
\end{simpro}

\begin{rem} Note that $\{\eta_i^{\ast}\}$ is a sequence of
independent random variables generated from a prespecified
distribution, and that its choice has little effect on the
finite-sample performance of the bootstrap simulation method. In
Example~\ref{ex4.1} we generate a sequence of $\{\eta_i^{\ast}\}$ from a
sequence of the normalized $\chi_2^2$ random variables of the form
$\frac{\chi_{2}^2 - 2}{2}$ even when $\{\eta_i\}$ is actually a
sequence of standard normal random variables. As discussed by Li and
Wang \cite{liwa98}, B\"{u}hlmann \cite{bu02}, Franke, Kreiss and Mammen \cite{frkrma02},
Chen and Gao \cite{chga07}, and others, in general we may use a wild bootstrap
method to generate a sequence of resamples for
$\{e_t^{\ast}\}$.
\end{rem}

As suggested by the referee, in Section~\ref{sec4.3} we propose using a block
bootstrap method (see, e.g., Hall, Horowitz and Jing \cite{hahoji95}) to compare
the finite-sample performance of the proposed test.

To investigate asymptotic properties of $l_{r}^{\ast}$ and
$\widehat{L}_{3n}^*(h)$, we need the following additional assumptions.
Assumption~\ref{ass4.1} ensures that
$\widehat{L}_{3n}(h)$ has some power under the alternative hypothesis.
Assumption~\ref{ass4.2} is a bootstrap version of Assumption~\ref{ass2.6}(i).
\begin{ass}\label{ass4.1}
 Let ${H}_0$ be false. \emph{(i)} Assumption~\ref{ass2.6}
holds, with $\theta_0$ replaced by
$\theta_1$.
\emph{(ii)} $\lim_{n\rightarrow\infty} n^{\alpha} h^{\alpha- 1/2} c_n^2=\infty$ for $\frac{1}{2}<\alpha<1$ and $0<\int_0^1
\Delta^2(x)
\,\mathrm{d}x<\infty$, where $c_n$ and $\Delta(x)$ are as defined in
(\ref{eq:1.2}).
\end{ass}
\begin{ass}\label{ass4.2} Under ${H}_0$,
$\Vert\widetilde{ \theta}^{\ast}-\widetilde{\theta}\Vert =\mathrm{O}_{P^{\ast
}}(n^{-1/2})$.
\end{ass}

We now have the following theorem.
Its proof is similar to that of Theorems~\ref{them2.1} and~\ref{them2.2} and is provided in
Appendix~\ref{appb}.

\begin{them}\label{them4.1} \emph{(i)} If, in addition to the conditions of
Theorem~\ref{them2.2},
Assumption~\ref{ass4.2} holds, then, under $H_{0}$,
%
\begin{equation}\label{eq:3.1.1}
\sup_{x \in R} \bigl|P^{\ast}\bigl(\widehat{L}_{3n}^*(h)\leq x\bigr) - P\bigl(\widehat
{L}_{3n}(h)\leq
x\bigr) \bigr| = \mathrm{o}_P(1),
\end{equation}
and $\lim_{n\rightarrow\infty} P (\widehat{L}_{3n}(h)>l_r^{\ast} )
= r$.

\emph{(ii)} If, in addition to the conditions of Theorem
\ref{them2.2},
Assumptions~\ref{ass4.1} and~\ref{ass4.2} hold, then, under ${H}_1$, $\lim_{n\rightarrow
\infty} P (\widehat{L}_{3n}(h)>l_r^{\ast} ) = 1$.
\end{them}

Theorem~\ref{them4.1} shows that the bootstrap approximation works well
asymptotically. For both the independent and short-range dependent
stationary cases, a discussion of existing results is available in
Chapter 3 of Gao \cite{ga07}. To the best of our knowledge, Theorem~\ref{them4.1} is a
new result in this kind of long-range dependent time series specification.

Note that $l_r^{\ast}$ is a function of $h$. A natural problem raised in
simulation is the choice of a suitable bandwidth $h$. To examine how to
address this issue, define the size and power functions of $\widehat
{L}_{3n}(h)$ as
%
\begin{equation}\label{eq:3.3}
\gamma_n(h) = P \bigl(\widehat{L}_{3n}^{\ast}(h) >l_r^{\ast}|{H}_0 \bigr)
\quad\mbox{and}\quad \beta_n(h) = P \bigl(\widehat{L}_{3n}^{\ast}(h)>l_r^{\ast
}|{H}_1 \bigr).
\end{equation}

Clearly, a reasonable selection procedure for a suitable bandwidth is such
that the size function $\gamma_n(h)$ is controlled by a significance level,
but the power function $\beta_n(h)$ is maximized over such bandwidths that
make $\gamma_n(h)$ is controllable. In this paper we propose choosing
%
\begin{equation}\label{eq:3.6}
\widehat{h}_{\mathrm{test}} = \arg\max_{h\in\mathcal{H}_n} \beta
_n(h) \qquad\mbox{with } \mathcal{H}_n = \{h\dvt r-\varepsilon<\gamma_n(h) < r
+ \varepsilon\},
\end{equation}
where $0<\varepsilon<r$. Theoretically, we have not been able to study
$\widehat{h}_{\mathrm{test}}$
asymptotically. In Example~\ref{ex4.1} below, we instead combine our proposed
Simulation Procedure~\ref{simpro4.1} and Implementation Procedure~\ref{imppro4.1} to numerically
approximate $\widehat{h}_{\mathrm{test}}$.
\begin{imppro}\label{imppro4.1} {Use $\widehat l_r^{\ast}=
l_r^{\ast}(\widehat{h}_{\mathrm{test}})$ as the simulated critical value
to compute the sizes and power values of $\widehat{L}_{3n}(\widehat
{h}_{\mathrm{test}})$}.
\end{imppro}
%
\subsection{LRD parameter estimation}\label{sec4.2} As mentioned at the beginning of
Section~\ref{sec4}, we need to estimate $\lambda=(\alpha,\eta)$ when $\lambda
$ is unknown. Let $u_t \equiv\widehat{e}_t =Y_t- m_{\widetilde
{\theta}}(X_t)$ and $I_u(\omega) = \frac{1}{2\uppi n} |\sum_{s=1}^n
u_s \mathrm{e}^{\mathrm{i}s \omega} |^2$.

Concerning the estimation of $\lambda$ for the case where the errors
are not necessarily Gaussian, several methods are available in the
literature. Giraitis and Surgailis \cite{gisu90} established an asymptotic
theory for a Whittle estimation method. Heyde and Gay \cite{hega93} considered
a multivariate case and established several asymptotic results.
Robinson \cite{ro95} proposed an efficient Gaussian semi-parametric
estimation method for the long-range dependence parameter, which is
equivalent to $\alpha$ in the present paper. Ideally, we would adopt
the estimation method proposed by Robinson \cite{ro95} for our case. But
because that method does not directly imply asymptotic normality for
$\widetilde{\eta}$ to meet our Assumption~\ref{ass2.3}, we modified the method
of Heyde and Gay \cite{hega93} for our case.

Define the following Gauss--Whittle objective function:
%
\begin{equation}
\Gamma_u(\lambda) = \frac{1}{4 \uppi} \int_{-\uppi}^{\uppi} \biggl(\log
(\psi_u(\omega;
\lambda)) + \frac{I_u(\omega)}{\psi_u(\omega; \lambda)} \biggr) \,\mathrm{d}\omega,
\label{eq:3.11a}
\end{equation}
where $\psi_u(\omega; \lambda)$ is the spectral density function of
$\{u_t\}$.

Let $\widetilde{\lambda}$ minimize $\Gamma_u(\lambda)$ over
$\Lambda_0$, a compact subset of $\Lambda=\{\lambda\dvt \lambda=
(\alpha, \eta)\dvt \frac{1}{2}<\alpha<1, 0<\eta<\infty\}$. Define
$L_i((-\uppi,\uppi]) = \{\psi\dvt \int_{-\uppi}^{\uppi} |\psi(\omega)|^i
\,\mathrm{d}\omega<\infty\}$ for $i=1,2$.

We introduce the following assumption:
\begin{ass}\label{ass4.3} \emph{(i)} Let $\{e_t\}$ satisfy Assumption~\ref{ass2.1} and
let $\psi_e(\omega, \lambda)$ be its spectral density function.

\emph{(ii)} Let $\lambda=(\lambda_1, \lambda_2)=(\alpha, \eta)$. The
functions $\psi_e(\omega; \lambda)$ and $p_i(\omega, \lambda) = -
\frac{\partial\psi_e^{-1}(\omega, \theta)}{\partial\lambda_i}$
for $i=1,2$ satisfy the following conditions:
\begin{itemize}
\item
$\int_{-\uppi}^{\uppi} \log(\psi_e(\omega, \lambda)) \,\mathrm{d}\omega$ is
twice differentiable in $\lambda$ under the integral sign. In
addition, $\psi_e(\omega, \lambda)$ is continuous at all $\omega
\neq0$ and $\lambda\in\Lambda$, and $\psi_e^{-1}(\omega, \lambda
)$ is continuous at all $(\omega, \lambda)$.
\item
$\psi_e^{-1}(\omega, \lambda)$, $\omega\in(-\uppi, \uppi]$ and
$\lambda\in\Lambda$, is twice differentiable with respect to
$\lambda$, and the functions $\frac{\partial\psi_e^{-1}(\omega,
\lambda)}{\partial\lambda_i}$ and $\frac{\partial^2 \psi
_e^{-1}(\omega, \lambda)}{\partial\lambda_i\, \partial\lambda_j}$
are continuous at all $(\omega, \lambda), \omega\neq0$ for all
$1\leq i, j\leq2$.
\item
For all $1\leq i\leq2$, $p_i(\omega, \lambda)$ are symmetric about
$\omega=0$ and $p_i(\omega, \lambda)\in L_1((-\uppi, \uppi])$ for any
$\lambda\in\Lambda$.
\item
$\psi_e(\omega, \lambda) p_i(\omega, \lambda)$ for $i=1,2$ are in
$L_2((-\uppi, \uppi])$ for any $\lambda\in\Lambda$. In addition, there
exists a constant $0<q\leq1$ such that $|\omega|^q \psi_e(\omega,
\lambda)$ is bounded and $|\omega|^{-q} p_i(\omega, \lambda)$ for
$i=1,2$ are in $L_2((-\uppi, \uppi])$ for any $\lambda\in\Lambda$.
\item
The matrix $ (\frac{\partial\log(\psi_e(\omega, \lambda
))}{\partial\lambda} ) (\frac{\partial\log(\psi_e(\omega,
\lambda))}{\partial\lambda} )^{\tau}$ is in $L_1((-\uppi, \uppi])$
for any $\lambda\in\Lambda$ and $\Sigma(\lambda) = \frac{1}{4 \uppi
} \int_{-\uppi}^{\uppi} (\frac{\partial\log(\psi_e(\omega, \lambda
))}{\partial\lambda} ) (\frac{\partial\log(\psi_e(\omega,
\lambda))}{\partial\lambda} )^{\tau} \,\mathrm{d}\omega$ is positive definite
for all $\lambda\in\Lambda$.
\end{itemize}
\end{ass}

Assumption~\ref{ass4.3} is equivalent to Conditions (A1)--(A3) of Heyde and Gay
\cite{hega93}. Assumption~\ref{ass4.3}(ii) is satisfied in many cases. Using existing
results (see, e.g., Beran \cite{be94}), we can verified that Assumption
\ref{ass4.3}(ii) is satisfied when the spectral density function is of the form
$\psi_e(\omega, \lambda) = d(\omega, \lambda) (\sin(\frac{\omega
}{2} ) )^{-(1-\alpha)}$, where $d(\omega, \lambda)>0$ satisfies
certain conditions.

Theorem~\ref{them4.2} establishes an asymptotic consistency result for $\tilde
\lambda$. Its proof is given in Appendix~\ref{appb}.

\begin{them}\label{them4.2}
Let the conditions of Theorem~\ref{them2.2}, except
Assumption~\ref{ass2.1}, hold. If in addition Assumption~\ref{ass4.3} holds, then, for
$n$ large enough,
%
\begin{equation}\label{eq:3.14}
\Vert\widetilde{\lambda} - \lambda\Vert =\mathrm{o}_P (\log n/n^{2/5} ),
\end{equation}
where $\widetilde{\lambda} = (\widetilde{\alpha}, \widetilde{\eta
} )$. Theorem~\ref{them4.2} shows that Assumption~\ref{ass2.3} can be justified for
$w_n=\frac{n^{2/5}}{\log(n)}$.
\end{them}
%
\subsection{An example of implementation}\label{sec4.3}
In this section we implement our proposed simulation procedure to show how
to assess the finite-sample properties of the proposed
test $\widehat{L}_{3n}(h)$ by using a simulated example.

Before we examine the finite-sample performance of the size and power
functions of our proposed test, we briefly state the following
simulation procedure as an alternative to the proposed Simulation
Procedure~\ref{simpro4.1}.
\begin{simpro}\label{simpro4.2}
Let $l_r$ ($0<r<1$) be the $1-r$
quantile of the exact finite-sample distribution of $\widehat
{L}_{3n}(h)$. Because $l_{r}$ might not be evaluated in practice, we
suggest an approximate $r$-level critical value, $l_{r}^{\ast}$, to
replace it, using the following bootstrap procedure:
\begin{enumerate}
\item
Generate $\{\widetilde{e}_t\}$ by $\widetilde{e}_t = \sum
_{j=1}^{\infty} \widetilde{\psi}_j \eta_{t-j}$, in which $\{\eta
_k\}$ is a sequence of independent observations drawn from $N(0,1)$,
and $\widetilde{\psi}_j=c(\widetilde{\alpha}) j^{-(1+
\widetilde{\alpha})/2}$, in which $c(\widetilde{\alpha}) = \sqrt
{\widetilde{\alpha}}$.
Let $l= [n^{1/3} ]$, and choose $b$ such that $b l =n$.
Generate $e_{1l}^{\ast}(j)= (\widetilde{e}_1(j), \ldots, \widetilde
{e}_l(j) ), \ldots$\,, and $e_{Nl}^{\ast}(j) = (\widetilde{e}_{(b-1)
l +1}(j), \ldots, \widetilde{e}_{bl}(j) )$ in step $j$ for $N = n - l
+1$. Replicate the resampling $J=250$ times and obtain $J$ bootstrap
resamples $\{e_{sl}^{\ast}(j)\dvt 1\leq s\leq N; 1\leq j\leq J\}$. Then
take the average $e_{sl}^{\ast} = \frac{1}{J} \sum_{j=1}^J
e_{sl}^{\ast}(j)$ to obtain a block bootstrap version of $e_t$ of the
form $ (e_1^{\ast}, \ldots, e_n^{\ast} ) = (e_{1l}^{\ast}, \ldots,
e_{Nl}^{\ast} )$.
Then generate $Y_t^{\ast}$ from $Y_t^{\ast}=m_{\widetilde{\theta
}}(X_t) + e_t^{\ast}$ for $1\leq t\leq n$, where $X_t=\frac{t}{n}$
for $1\leq t\leq n$.
\item
Use the data set $\{Y_t^{\ast}\dvt 1\leq t\leq n\}$ to obtain an
estimator $\widetilde{\theta}^{\ast}$ of ${\widetilde\theta}$, and
construct a corresponding version $\widehat{L}_{1n}^*(h)$ of $\widehat
{L}_{1n}(h)$ under $H_0$; that is, compute $\widehat{L}_{1n}^*(h)$
according to the same formula as for $\widehat{L}_{1n}(h)$, but with
$\{Y_t\dvt 1\leq t\leq n\}$ and ${(\theta_0,\widetilde\theta)}$
replaced by $\{Y_t^{\ast}\dvt 1\leq t\leq n\}$ and $(\widetilde\theta
,\widetilde{\theta}^{\ast})$, respectively.
\item
Repeat the foregoing step $M$ times and produce $M$ versions of
$\widehat{L}_{3n}^{\ast}(h)$, denoted by $\widehat{L}_{3n,m}^{\ast
}(h)$ for $m=1,2,\ldots, M$. Use the $M$ values of $\widehat
{L}_{1n,m}^{\ast}(h)$ to construct the empirical distribution
function. Define $P^{\ast}(\cdot)$ as the bootstrap distribution of
$\widehat{L}_{1n}^*(h)$ given\vspace*{1pt} $\mathcal{W}_n$ by $P^{\ast} (\widehat
{L}_{1n}^*(h)\leq x )=P (\widehat{L}_{1n}^*(h)\leq x|\mathcal{W}_n
)$. Let $l_{r}^{\ast}$ ($0<r<1$) satisfy $P^{\ast} (\widehat
{L}_{1n}^*(h)\geq
l_{r}^{\ast} )=r$, and estimate $l_{r}$ by $l_{r}^{\ast}$.
\end{enumerate}
\end{simpro}

\begin{ex}\label{ex4.1} Consider a linear model of the form
%
\begin{equation}\label{eq:4.1}
Y_t = \alpha_0 + \beta_0 X_t + e_t,\qquad t=1, \ldots, n,
\end{equation}
where $(\alpha_0, \beta_0)$ is a pair of unknown parameters,
$X_t =\frac{t}{n}$ for $1\leq t\leq n$, and $\{e_t\}$ is a sequence of
dependent errors given by
$e_t = \sum_{j=1}^{\infty} \psi_j \eta_{t-j}$,
where $\{\eta_k\}$ is a sequence of independent
observations drawn from $N(0,1)$ and $\psi_j=c(\alpha) j^{-(1+\alpha)/2}$ for $\frac{1}{2}<\alpha<1$, with $c(\alpha)=\sqrt
{\alpha}$. We choose $\alpha= \frac{3}{4}$ in our simulation.
Throughout this section, we use the standard normal kernel function
$K(x) = \frac{1}{\sqrt{2\uppi}} \mathrm{e}^{-x^2/2}$.

To compute the sizes and power values of $\widehat{L}_{3n}(h)$, we
generate $\{Y_t\}$ from
%
\begin{equation}
H_0\dvt Y_t = \alpha_0 + \beta_0 X_t + e_t \quad\mbox{or}\quad H_1\dvt Y_t =
\alpha_1 + \beta_1 X_t + \gamma_1 X_t^2 + e_t,
\label{eq:5.4}
\end{equation}
where the parameters $(\alpha_0, \beta_0)$ are estimated by
$(\widetilde{\alpha}_0, \widetilde{\beta}_0)$ under $H_0$
and the parameters $(\alpha_1, \beta_1, \gamma_1)$ are estimated
by the ordinary least squares estimators $(\widetilde{\alpha}_1,
\widetilde{\beta}_1, \widetilde{\gamma}_1)$ under $H_1$. When we generate
$\{Y_t\}$, the initial values are $\alpha_i=\beta_i\equiv1$ for
$i=0,1$, and we use $\gamma_1= n^{-\frac{1}{2}} \sqrt{\mathrm
{loglog(\mathit{n})}}$ to compute the power values in both Tables \hyperref[table1]{1} and \hyperref[table2]{2}.
\begin{table}
\tablewidth=\textwidth
\tabcolsep=0pt
  \caption{Sizes and power values based on simulated critical values}\label{table1}
\begin{tabular*}{\textwidth}{@{\extracolsep{\fill}}llllllll@{}}
 \multicolumn{8}{c}{Null hypothesis is true}
 \\
 \hline
 && \multicolumn{2}{l}{$r=1\%$} & \multicolumn
{2}{l}{$r=5\%$} & \multicolumn{2}{l}{$r=10\%$}
\\[-5pt]
&& \multicolumn{2}{l}{\hrulefill} & \multicolumn
{2}{l}{\hrulefill} & \multicolumn{2}{l@{}}{\hrulefill}
\\
{Observation}&$T$: & {$f_{0\mathrm{cv}}$} & {$f_{0\mathrm{test}}$} & {$f_{0\mathrm
{cv}}$} & {$f_{0\mathrm{test}}$} & {$f_{0\mathrm{cv}}$} &
{$f_{0\mathrm{test}}$}
\\
\hline
{250} & & {0.004} & 0.011 & {0.026} & 0.064 & {0.061} & 0.086
\\
{500} & & {0.006} & 0.009 & {0.031} & 0.060 & {0.058} & 0.082
\\
{750} & & {0.003} & 0.012 & {0.028} & 0.042 & {0.049} & 0.122
\\
\hline  \\[-2pt]
\multicolumn{8}{c}{Null hypothesis is false}
\\
\hline
& & \multicolumn{2}{l}{$r=1\%$} & \multicolumn
{2}{l}{$r=5\%$} & \multicolumn{2}{l}{$r=10\%$}
\\[-5pt]
& & \multicolumn{2}{l}{\hrulefill} & \multicolumn
{2}{l}{\hrulefill} & \multicolumn{2}{l@{}}{\hrulefill}
\\
{Observation} & $T$: & {$f_{1\mathrm{cv}}$} & {$f_{1\mathrm{test}}$} & {$f_{1\mathrm
{cv}}$} & {$f_{1\mathrm{test}}$} & {$f_{1\mathrm{cv}}$} &
{$f_{1\mathrm{test}}$}
\\
\hline
{250} & & 0.084 & 0.160 & 0.121 & 0.242 & 0.181 & 0.366
\\
{500} & & 0.076 & 0.144 & 0.153 & 0.296 & 0.179 & 0.304
\\
{750} & & 0.104 & 0.216 & 0.142 & 0.320 & 0.234 & 0.426
\\
\hline
\end{tabular*}
\end{table}

\begin{table}[b]
\tablewidth=\textwidth
\tabcolsep=0pt
  \caption{Sizes and power values based on simulated critical values}\label{table2}
\begin{tabular*}{\textwidth}{@{\extracolsep{\fill}}llllllll@{}}
\multicolumn{8}{c}{Null hypothesis is true}
\\
\hline
& & \multicolumn{2}{l}{$r=1\%$} & \multicolumn
{2}{l}{$r=5\%$} & \multicolumn{2}{l}{$r=10\%$}
\\[-5pt]
& & \multicolumn{2}{l}{\hrulefill} & \multicolumn
{2}{l}{\hrulefill} & \multicolumn{2}{l@{}}{\hrulefill}
\\
{Observation}& $T$: & {$f_{0\mathrm{cv}}$} & {$f_{0\mathrm{test}}$} & {$f_{0\mathrm
{cv}}$} & {$f_{0\mathrm{test}}$} & {$f_{0\mathrm{cv}}$} &
{$f_{0\mathrm{test}}$}
\\
\hline
{250} & & {0.004} & 0.010 & {0.016} & 0.060 & {0.050} & 0.136
\\
{500} & & {0.002} & 0.012 & {0.016} & 0.062 & {0.034} & 0.113
\\
{750} & & {0.006} & 0.010 & {0.018} & 0.051 & {0.036} & 0.102
\\
\hline\\[-2pt]
\multicolumn{8}{c}{Null hypothesis is false}
\\
\hline
& & \multicolumn{2}{l}{$r=1\%$} & \multicolumn
{2}{l}{$r=5\%$} & \multicolumn{2}{l}{$r=10\%$}
\\[-5pt]
& & \multicolumn{2}{l}{\hrulefill} & \multicolumn
{2}{l}{\hrulefill} & \multicolumn{2}{l@{}}{\hrulefill}
\\
{Observation}& $T$: & {$f_{1\mathrm{cv}}$} & {$f_{1\mathrm{test}}$} & {$f_{1\mathrm
{cv}}$} & {$f_{1\mathrm{test}}$} & {$f_{1\mathrm{cv}}$} &
{$f_{1\mathrm{test}}$}
\\
\hline
{250} & & 0.080 & 0.168 & 0.130 & 0.324 & 0.190 & 0.422
\\
{500} & & 0.070 & 0.332 & 0.164 & 0.526 & 0.224 & 0.628
\\
{750} & & 0.114 & 0.334 & 0.180 & 0.544 & 0.246 & 0.654
\\
\hline
\end{tabular*}
\end{table}

We first apply Implementation Procedure~\ref{imppro4.1} for the case where $\{\eta
_i^{\ast}\}$ is a sequence of the normalized $\chi_2^2$ random
variable of the form $\frac{\chi_2^2 -2}{2}$. When using
Simulation Procedure~\ref{simpro4.2}, we start with the case where $\{\eta_i\}
$ is a sequence of i.i.d. observations drawn from $N(0,1)$. We then
apply either Implementation Procedure~\ref{imppro4.1} or Simulation
Procedure~\ref{simpro4.2} to find $\widehat{h}_{\mathrm{test}}$.\vadjust{\goodbreak} To assess whether the use
of $\widehat{h}_\mathrm{test}$ associated with the proposed bootstrap
method could improve the finite-sample performance of $\widehat
{L}_{3n}$, we also consider using an estimation-based optimal bandwidth
of the form $\widehat{h}_{{ cv}} = \arg\min_{h\in H_\mathrm{cv}}
\frac{1}{n} \sum_{t=1}^n (Y_t - \widehat{m}_{-t}(X_t,h) )^2$,
where $\widehat{m}_{-t}(X_t,h) = \sum_{s=1, \neq t}^n K (\frac
{X_t - X_s}{h} )
Y_s/\sum_{u=1, \neq t}^n K (\frac{X_t - X_u}{h} )$ and $H_\mathrm
{cv} = [c_1 n^{-1}, c_2 n^{-(1-c_0)} ]$ for some $0<c_1<c_2<\infty$
and $0<c_0<1$.

In the implementation of Simulation Procedure~\ref{simpro4.1}, we consider the case
with $M = 500$ replications, each with $B=250$ bootstrapping resamples.
When implementing Simulation Procedure~\ref{simpro4.2}, we use the same number of
replications $M=500$, each with $J=250$ block bootstrapping resamples
($n = bl$). All of these simulations were done for data sets of size
$n=250$, $500$ and~$750$. Let $\widehat{h}_{i\mathrm{test}}$ and
$\widehat{h}_{i\mathrm{cv}}$ denote the corresponding versions of
$\widehat{h}_\mathrm{test}$ and $\widehat{h}_\mathrm{cv}$ under
$H_i$ for
$i=0,1$. Let $f_{i\mathrm{test}}$ denote the frequency of
$\widehat{L}_{3n}(\widehat{h}_{i\mathrm{test}})>l_{r}^{\ast
}(\widehat{h}_{i\mathrm{test}})$ (note that we repeat the test $10$
times to compute $l_{r}^{\ast}(\widehat{h}_{i\mathrm{test}})$), and let
$f_{i\mathrm{cv}}$ be the frequency of
$\widehat{L}_{3n}(\widehat{h}_{i\mathrm{cv}})>z_r$ under $H_{i}$ for each
of $i=0,1$, where $z_{0.01}=2.33$ at the $1\%$ level,
$z_{0.05}=1.645$ at the $5\%$ level, and $z_{0.10}=1.28$ at the
$10\%$ level.

Tables \hyperref[table1]{1} and \hyperref[table2]{2} present the simulated sizes and power values based
on Simulation Procedures~\ref{simpro4.1} and~\ref{simpro4.2}, respectively.

Tables \hyperref[table1]{1} and \hyperref[table2]{2} show that the simulated power values associated with
$\widehat{h}_\mathrm{test}$ are greater than those based on $\widehat
{h}_\mathrm{cv}$, and also that $\widehat{L}_{1n}(\widehat
{h}_\mathrm{cv})$ has some kind of size distortion when using $z_r$
(equivalent to using the asymptotic normality as the sample
distribution) in practice. This finding is not surprising, given that
the the theory demonstrates that each $\widehat{h}_{\mathrm{test}}$
is chosen such to maximize the resulting power function, with the
corresponding size function computed using the simulated critical value
$l_r^{\ast}(\widehat{h}_\mathrm{test})$ in each case.

Comparing the regression bootstrap method in Table \hyperref[table1]{1} and the block
bootstrap method in Table~\hyperref[table2]{2} shows similar sizes for the two methods.
Meanwhile, the power values in columns 3, 5 and~7 in Tables~\hyperref[table1]{1}
and~\hyperref[table2]{2} show that $\widehat{L}_{3n}(\widehat{h}_\mathrm{test})$
associated with the block-bootstrap method is more powerful than
$\widehat{L}_{3n}(\widehat{h}_\mathrm{test})$ based on the
regression bootstrap method, whereas whether $\widehat
{L}_{3n}(\widehat{h}_\mathrm{cv})$ associated with the
block-bootstrap method or $\widehat{L}_{3n}(\widehat{h}_\mathrm
{cv})$ based on the regression bootstrap method is uniformly more
powerful is not so clear. This finding further supports our argument
that the choice of a suitable bootstrap method for selecting an
appropriate bandwidth to maximize the resulting power function is more
relevant than the choice of the bootstrap method itself.
\end{ex}

\section{Conclusion}\label{sec5}
 We have proposed a new nonparametric test for the
parametric specification of the mean function of long-range dependent time
series, and have established asymptotic distributions of the proposed
test for both the fixed and
random design cases. In addition, we have proposed Simulation Procedure
\ref{simpro4.1} and Implementation Procedure
\ref{imppro4.1} to implement the proposed test in practice. Our finite-sample results
show that the proposed test, as well as the two procedures, are practically
applicable and implementable. Further topics, including how to extend
existing results (Nishiyama and Robinson \cite{niro00,niro05}; Gao \cite{ga07}; Gao and
Gijbels \cite{gagi08}) for Edgeworth expansions for the size and power functions
of the proposed tests to the long-range dependence case, are left for
future research.
\begin{appendix}
\section{Technical details}\label{appa}
This appendix provides technical details for the asymptotic theory in
Section~\ref{sec2}. Appendix~\ref{appa.1} gives several preliminary lemmas. Appendix~\ref{appa.2}
presents the proofs of Theorems~\ref{them2.3} and~\ref{them2.4}, and~\ref{appa.3} presents the
proofs of Theorems~\ref{them2.1} and~\ref{them2.2}. Proofs of the lemmas are given in
Appendix~A of Gao and Wang \cite{gawa10}. Throughout this section, we denote
constants by $C, C_1,\ldots$\,, which may take different values at each
appearance.
%
\subsection{Technical lemmas}\label{appa.1}
\begin{lem}\label{lema.1}
Let $\{e_t\}$ be a linear process defined
by $e_t=\sum_{j=-\infty}^{\infty}\psi_j\eta_{t-j},
$ where $\{\eta_j\}$ is a sequence of i.i.d. random variables with
$E[\eta_1]=0$,
$E[\eta_1^2]=1$,
$E[\eta_1^4]<\infty$ and $\gamma(0)<\infty$, where
$\gamma(k)=E[e_t e_{k+t}]=\sum_{j=-\infty}^{\infty}\psi_j\psi_{j+k}$.
Then, for all $j, k, s$ and $t$,
%
\begin{eqnarray}\label{ad10}
E[e_j e_k e_s e_t]
&=& (E[\eta_1^4] - 3) \sum_{m=-\infty}^{\infty}
\psi_{j-m}\psi_{k-m}\psi_{s-m}\psi_{t-m}\no
\\[-8pt]
\\[-8pt]
&&{} + \gamma(j-k)\gamma(s-t)+\gamma(j-s)\gamma(k-t)+\gamma
(j-t)\gamma(k-s).\no
\end{eqnarray}
In particular, we have that $E [e_1^4 ]<\infty$,
%
\begin{eqnarray}
\label{fin15}|E [e_j^2e_k^2 ]- \gamma^2(0)| &\le& C \gamma^2(j-k) \quad\mbox{and}\\
 \label{fin1}|E [e_j^2 e_k e_s ] |
&\leq& C [\gamma(j-k)+\gamma(j-s)+\gamma(k-s) ]
\end{eqnarray}
for all $j\not= k\not= s$. If, in addition, $\psi_k\ge0$,
then for all $j, k, s$ and $t$,
%
\begin{eqnarray}\label{fin1a9}
&& \bigl|E \bigl[ \bigl(e_j e_{j+s}-\gamma(s) \bigr) \bigl(e_k
e_{k+t}-\gamma(t)\bigr )\bigr ]
\bigr| \no
\\[-8pt]
\\[-8pt]
&& \quad\le C \gamma(j-k)\gamma(j-k+s-t)+\gamma(j-k+s)\gamma(j-k-t).\no
\end{eqnarray}
\end{lem}
\begin{lem}\label{lema.2} Let $1/2<\al<1$ and $0<\beta< \al-1/2$.
Then, for all
$k\ge3$ and as $n\to\infty$,
%
\begin{eqnarray}\label{fin1a}
I_n&=& \f1{n^{k/2}}\int_{1}^n \int_1^{n}\cdots\int_1^n
|x_1-x_2|^{-\al}|x_2-x_3|^{\beta-1} \cdots\nonumber
\\[-8pt]
\\[-8pt]
&&\hphantom{\f1{n^{k/2}}\int_{1}^n \int_1^{n}\cdots\int_1^n}{} \times |x_{2k-1}-x_{2k}|^{-\al}
|x_{2k}-x_1|^{\beta-1} \,\mathrm{d}x_1\,\mathrm{d}x_2\cdots \,\mathrm{d}x_{2k} \to0.\nonumber
\end{eqnarray}
\end{lem}

In Lemma~\ref{lema.3} below, let $\{X_i, i\ge1\}$ be a sequence of i.i.d. random
variables with density function $f(x)$, and set $g_n(X_i, X_j)=K (\frac
{X_i - X_j}{h} ) - E [K (\frac{X_i - X_j}{h} ) ]$,
%
\begin{eqnarray}
\label{hk1}
g_{1n}(X_i) & = & E [g_n(X_i, X_j)|X_i ] \quad\mbox{and} \nonumber
\\[-8pt]
\\[-8pt]
g_{2n}(X_i, X_j) & = & g_n(X_i, X_j) -g_{1n}(X_i)-g_{1n}(X_j).\nonumber
\end{eqnarray}

\begin{lem}\label{lema.3}
 Let $K(x)$ satisfy Assumption~\ref{ass2.4}\emph{(i)}.
If $f(x)$ is a bounded and uniformly continuous function on $\R$, then
%
\begin{eqnarray}
\label{try1}E \bigl[K [(X_1- X_2)/h ] \bigr] &\sim& c_1 h,  \\
\label{try2}E [g_{2n}^2(X_1,X_2) ] &\sim& E \bigl[K^2 [(X_1-X_2)/h ]
\bigr] \sim c_2 h, \\
\label{try3}E [g_{2n}^4(X_1,X_2) ] &\sim&
E \bigl[K^4 [(X_1-X_2)/h ] \bigr] \sim c_4 h,
\end{eqnarray}
where $c_j=\int_{-\infty}^{\infty} K^j(s)\,\mathrm{d}s \int_{-\infty}^{\infty}
f^2(y)\,\mathrm{d}y$ for $j=1, 2, 4$.
Furthermore,
%
\begin{eqnarray}
\label{hk4}E [g_{1n}^2(X_1) ] & \sim& d_1 h^2,  \\
\label{hk2}E [g_{2n}(X_1,X_3) g_{2n}(X_1,X_{4}) g_{2n}(X_2,X_3) g_{2n}(X_2,X_{4}) ]
& \sim&
d_2 h^3,
\end{eqnarray}
where
\[
d_1 = \int_{-\infty}^{\infty}
\biggl[f(x)-\int_{-\infty}^{\infty}f^2(y)\,\mathrm{d}y \biggr]^2f(x)\,\mathrm{d} x
\]
and
\[
d_2 =\int\int\int K(s)K(t)
K(x+s)K(x+t)\,\mathrm{d}s\,\mathrm{d}t\,\mathrm{d}x \int_{-\infty}^{\infty}f^4(y)\,\mathrm{d}y.
\]
\end{lem}

Our next lemma establishes a Berry--Esseen type bound for random weighted
$U$ statistics.
This lemma is interesting and useful in itself.
%
\begin{lem}\label{lema.4}
Let $\{\varepsilon_k, k\ge1\}$ be a sequence
of i.i.d. random variables, let $\{a_{nij}\}$ be a sequence of constants
with $a_{nij}=a_{nji}$ for all $n\ge1$, and let $\{\varphi_n(x, y)\}$
be a sequence of symmetric Borel-measurable functions such that for all
$n\ge1$,
%
\begin{equation}\label{jan7}
E [\varphi_n^2(\ep_1,\ep_2) ]>0,\qquad
E [\varphi_n(\ep_1,\ep_2)\vert\ep_1 ]=0.
\end{equation}
Then
there exists an absolute constant $A>0$ such that
%
\begin{equation} \label{in2}
\sup_x |P (B_n^{-1}S_n\le x )-\Phi(x) |
\le A B_n^{-4/5} \bigl(A_{1n} E\varphi_n^4(\ep_1,\ep_2)
+A_{2n} \mathcal{L}_n \bigr)^{1/5},
\end{equation}
where $S_n
=\sum_{1\le i<j\le n} a_{nij} \varphi_n(\ep_i,\ep_j), $
$B_n^2=\sum_{1\le i<j\le n} a_{nij}^2E\varphi_n^2(\ep_1,\ep_2)$,
%
\begin{eqnarray}
A_{1n}&=&\sum_{i=2}^n \Biggl(\sum_{j=1}^{i-1}a_{nij}^2 \Biggr)^2,\qquad
A_{2n} = \sum_{i=2}^{n-1}\sum_{j=i+1}^{n}\Biggl (\sum_{k=1}^{i-1} a_{nik}
a_{njk} \Biggr)^2,
\\
\mathcal{L}_n &=& E [\varphi_n(\ep_1,\ep_3)
\varphi_n(\ep_1,\ep_{4}) \varphi_n(\ep_2,\ep_3) \varphi_n(\ep
_2,\ep_{4}) ].
\end{eqnarray}
\end{lem}
\begin{lem}\label{lema.5}
Let $\{\varepsilon_k, k\ge1\}$ be a sequence of
i.i.d. random variables with $E[\ep_1]=0$ and $E[\ep_1^6]<\infty$. Let
$\{a_{nij}\}$ be a sequence of real numbers
with $a_{nij}=a_{nji}$ and $\Vert A\Vert^2\equiv
\sum_{i,j=-\infty}^{\infty}a_{nij}^2<\infty$
for all $n\ge1$. If there exists an absolute constant $b_1^2>0$
such that $1-\frac{V^2}{\Vert A\Vert^2}\ge b_1^2$ with
$V^2=\sum_{i=-\infty}^{\infty}a_{nii}^2$, then
%
\begin{eqnarray}
\sup_x |P (S_n/B_n\le x )-\Phi(x) |
&\le& C \frac{ \{\operatorname{Tr}(A^4) \}^{1/4}}{\Vert A\Vert}, \label{in2aa}
\end{eqnarray}
where $S_n
=\sum_{i,j=-\infty}^{\infty} a_{nij} (\ep_i \ep_j-E[\ep_i\ep_j] ),
$
$B_n^2=2(\Vert A\Vert^2-V^2)\mu_2^2+V^2(\mu_4-\mu_2^2)$ with $\mu_j=E
[|\ep_1|^j ]$,
and $A$ is the infinite matrix with $a_{nij}$ as its $(i,j)$th element.
\end{lem}
\begin{lem}\label{lema.6}
Let $K(x)$ be a non-negative symmetric integrable
function
satisfying $K(x) = \mathrm{O} [(1+|x|^{1-\beta})^{-1} ]$,
where $0<\beta\le\al-1/2$ and $1/2<\al<1$. Then
%
\begin{eqnarray}
\label{last0}\Delta_0 &\equiv& \int_0^{\infty}x^{-\al} K(x)
\,\mathrm{d}x<\infty, \\
\label{last}A_{\al} &\equiv& \int_{0}^{\infty} \int_{0}^{\infty}
\int_{0}^{\infty} x^{-\al} y^{-\al} [I_1(x,y,w)+I_2(x,y,w) ]
\,\mathrm{d}x\,\mathrm{d}y\,\mathrm{d}w
<\infty,
\end{eqnarray}
and as $h\to0$,
%
\begin{eqnarray}\label{last1}
\Delta_1 &\equiv& \int_{0}^{1/h} \int_{0}^{1/h}
\int_{x}^{1/h} x^{-\al} y^{-\al} \max\{w, y\} [I_1(x,y,w) + I_2(x,y,w) ]
 \,\mathrm{d}x\,\mathrm{d}y\,\mathrm{d}w\no
 \\[-8pt]
\\[-8pt]
& = &\mathrm{o}(1/h),\no
\end{eqnarray}
where $I_1(x,y,w)=K(w)K(x+y-w)$ and $I_2(x,y,w)=K(w-x)K(w-y)$.
\end{lem}
%
\subsection{\texorpdfstring{Proofs of Theorems \protect\ref{them2.3} and \protect\ref{them2.4}}
{Proofs of Theorems 2.3 and 2.4}}\label{appa.2}
\begin{pf*}{Proof of Theorem~\ref{them2.3}}
 We may write
\[
\sum_{1\le i\not=j\le n} e_i e_j
\biggl[K \biggl(\frac{(X_i- X_j)}{h} \biggr) \biggr] = \widetilde
Q_{2n}^{(1)}+\widetilde Q_{2n}^{(2)}+\widetilde Q_{2n}^{(3)},
\nonumber
\]
where
\begin{eqnarray}
\widetilde Q_{2n}^{(1)} &=&\sum_{1\le i\not=j\le n} e_i e_j
E \biggl[K \biggl(\frac{(X_i- X_j)}{h}\biggr ) \biggr], \\
\widetilde Q_{2n}^{(2)} &=&\sum_{1\le i\not=j\le n} e_i e_j [
g_{1n}(X_i)+g_{1n}(X_j) ], \\
\widetilde Q_{2n}^{(3)} &=&\sum_{1\le i\not=j\le n} e_i e_j g_{2n}(X_i,
X_j),
\end{eqnarray}
where $g_{1n}(X_i)$ and
$g_{2n}(X_i, X_j)$ are defined as in (\ref{hk1}). Theorem~\ref{them2.3}
now follows readily if we prove the following: Whenever $h\to0$,
%
\begin{eqnarray}
(A_{2\al}n^{2-\al} h )^{-1} {\widetilde Q_{2n}^{(1)}} &\rightarrow_D&
\chi^2(1), \label{rel21} \\
\widetilde Q_{2n}^{(2)} &=& \mathrm{o}_P \bigl(\max\bigl\{n^{2-\al} h, n\sqrt
h \bigr\} \bigr), \label{rel22}
\end{eqnarray}
and if in addition, $nh\to\infty$, then
%
\begin{equation}
\bigl( A_{1\al} n \sqrt{h} \bigr)^{-1}
{\widetilde Q_{2n}^{(3)}} \rightarrow_D N(0,1). \label{rel23}
\end{equation}
%

Actually, if $h\to0$ and $\sqrt h n^{1-\al}\to\infty$, then
$\widetilde
Q_{2n}^{(2)} +\widetilde Q_{2n}^{(3)}=\mathrm{o}_P (n^{2-\al}h )$ by virtue of
(\ref{rel22}) and (\ref{rel23}).
This and (\ref{rel21}) yield
Theorem~\ref{them2.3}(ii). Similarly, if
$\sqrt h n^{1-\al}\to0$ and $nh\to\infty$, then
$\widetilde Q_{2n}^{(1)} +\widetilde Q_{2n}^{(2)}=\mathrm{o}_P (n\sqrt h )$ by
virtue of (\ref{rel21}) and (\ref{rel22}).
This and (\ref{rel23}), yield
Theorem~\ref{them2.3}(i).

We now prove (\ref{rel21})--(\ref{rel23}), starting with (\ref
{rel21}). By (\ref
{try1}),
\[
(A_{2\al}n^{2-\al} h )^{-1} {\widetilde Q_{2n}^{(1)}}
= \bigl(1+\mathrm{o}_P(1)\bigr) \Biggl[ \Biggl(\frac1{d_n}\sum_{j=1}^n e_j \Biggr)^2-\frac
1{d_n^2}\sum_{j=1}^n e_j^2 \Biggr],\label{hy8}
\]
where $d_n^2=\f
{2\eta}{(1-\al)(2-\al)}n^{2-\al}$.
It can be readily seen that $\frac1{d_n^2}\sum_{j=1}^n e_j^2\to0$
a.s. by the
stationary ergodic theorem. This, together with (\ref{hy8}) and the
continuous mapping theorem, yield that (\ref{rel21}) will follow
if we prove
%
\begin{equation}\label{ad17}
\frac1{d_n}\sum_{j=1}^n e_j\rightarrow_D N(0,1), \label{ad17}
\end{equation}
which follows from Lemma 1 of Robinson \cite{ro97}. This
proves (\ref{rel21}).

We next prove (\ref{rel22}). By (\ref{hk4}), independence of
$e_i$ and $X_i$, and (\ref{fin1}),
\begin{eqnarray*}\label{re5}
E \bigl(\widetilde Q_{2n}^{(2)}\bigr )^2 &=& E \Biggl(\sum_{i=1}^ng_{1n}(X_i)
e_i\sum_{1\le j
\not= i\le n}e_j\Biggr )^2\no\\
&=& [d_1 h^2+\mathrm{o}(h^2)] \biggl(\sum_{1\le i\not=j\le
n} E [e_i^2 e_j^2 ] +
\sum_{1\le i\not=k\not= j\le n}E [e_i^2 e_j e_k ] \biggr)\no\\
&\le& C h^2 \Biggl[n^2 E e_1^4+ n\sum_{i,j=1}^n\gamma(i-j) \Biggr] \le
C h^2 (n^2+n^{3-\al} ).
\end{eqnarray*}
Thus, equation (\ref{rel22}) follows immediately from the Markov
inequality.

Finally, we prove (\ref{rel23}). Write
\[
B_n^2=\sum_{1\le i<j\le n} e_i^2 e_j^2 E [g_{2n}^2(X_i,
X_j) ]
=\f12 E [g_{2n}^2(X_1, X_2) ] \cdot\Biggl[ \Biggl(\sum_{i=1}^n
e_i^2 \Biggr)^2-\sum_{i=1}^n e_i^4 \Biggr].
\]
By (\ref{try2}) and the stationary ergodic theorem, which yields that
$\frac{1}{n} \sum_{i=1}^n e_i^2\to E[e_1^2]=\gamma(0)$ a.s. and
$\frac{1}{n} \sum_{i=1}^n e_i^4 \to E[e_1^4]<\infty$ a.s., it can be readily
seen that as $n\to\infty$
%
\begin{equation}\label{lk1}
4 A_{1\al}^{-2} n^{-2} h^{-1}
{B_n^2}\to
1\qquad \mbox{a.s.},
\end{equation}
where $A_{1\al}$ is defined as in Theorem 2.3.
Thus, to prove (\ref{rel23}),
it suffices to show that
%
\begin{equation}\label{lk2}
\widetilde Q_{2n}^{(3)}/ (2 B_n) \rightarrow_D N(0,1).
\end{equation}

Lemma~\ref{lema.4} is used to establish (\ref{lk2}).
In fact, noting that
\[
\widetilde Q_{2n}^{(3)}=2 \sum_{1\le i<j\le
n} e_i e_j g_{2n}(X_i, X_j)
\]
and
$E [g_{2n}(X_i,X_j)|e_i ]=0$ for all $i\not=j$, it
follows from the independence of $e_i$ and $X_i$,
Lemma~\ref{lema.4}, (\ref{try3}) and (\ref{hk2}) that
%
\begin{equation}\label{in2a}
\sup_x \bigl|P \bigl(\widetilde Q_{2n}^{(3)}/2B_n\le x|e_1,
\ldots,e_n \bigr)-\Phi(x) \bigr|
\le A B_n^{-4/5} (c_4 h A_{1n}+d_2 h^3 A_{2n} )^{1/5},
\end{equation}
where $A$ is an absolute constant, $c_4$ and $d_2$ are defined as in
(\ref
{try3})
and (\ref{hk2}), and
\begin{eqnarray*}
A_{1n}&=&\sum_{i=2}^n
\Biggl(\sum_{j=1}^{i-1}(e_i e_j)^2 \Biggr)^2 \le\sum_{i=2}^n e_i^4\Biggl (\sum
_{j=1}^{n} e_j^2 \Biggr)^2,\\
A_{2n} &=&
\sum_{i=2}^{n-1}\sum_{j=i+1}^{n} \Biggl(\sum_{k=1}^{i-1} e_i e_j e_k^2\Biggr )^2
\le\biggl(\sum_{1\le i<j\le n} e_i^2 e_j^2 \biggr)^2.
\end{eqnarray*}
Again by the stationary ergodic theorem, for $n$ large enough,
\[
\frac{1}{n^3} A_{1n} \le2 E[e_1^4] \cdot(E[e_1^2])^2\qquad\mbox{a.s.}\quad
\mbox{and}\quad
\frac{1}{n^4} A_{2n} \le2 (E[e_1^2] )^4\qquad\mbox{a.s.}
\]
This, together with (\ref{lk1}) and (\ref{in2a}),
implies that for $n$ large enough,
\[
\sup_x \bigl|P\bigl (\widetilde Q_{2n}^{(3)}/(2B_n)\le
x|e_1,\ldots,e_n \bigr)-\Phi(x) \bigr|\le C
\biggl(\frac{1}{nh}+h \biggr)^{1/5}\qquad
\mbox{a.s.}
\]
Now if $h\to0$ and
$nh\to\infty$, then
\begin{eqnarray*}
&&\lim_{n\to\infty}\sup_x \bigl|P \bigl(\widetilde Q_{2n}^{(3)}/(2B_n)\le
x \bigr)-\Phi(x) \bigr| \\
&&\quad \le E \Bigl[\lim_{n\to\infty}\sup_x \bigl|P \bigl(\widetilde
Q_{2n}^{(3)}/(2B_n)\le
x\vert e_1, \ldots,e_n\bigr )-\Phi(x) \bigr| \Bigr]
=0.
\end{eqnarray*}
This proves (\ref{lk2}), and also completes the proof of
Theorem 2.3.
\end{pf*}
\begin{pf*}{Proof of Theorem 2.4}
 Let
\[
\widetilde Q_{1n}=K(0)\sum_{i=1}^n e_i^2 \quad\mbox{and}\quad
Q_{1n}= \sum_{i=1}^n\sum_{j=1, j\not=i}^n e_i e_j b_n(i,j).
\]
We
have that
%
\begin{eqnarray}
\widetilde{Q}_{n} &\equiv& Q_{1n} - E[Q_{1n}] +\widetilde Q_{1n} -
E[\widetilde{Q}_{1n}]
= \sum_{i,j=1}^n \bigl(e_i e_j
-\gamma(i-j) \bigr) b_n(i,j)\no
\\[-8pt]
\\[-8pt]
&=&\sum_{k,l=-\infty}^{\infty} a_{nkl} (\eta_k\eta_l-E [\eta
_k\eta_l ] ),\no
\end{eqnarray}
where $a_{nkl}=\sum_{i,j=1}^n\psi_{i-k}\psi_{j-l} b_n(i,j)$, based
on the fact that $E[\eta_k\eta_l]=0$ for $k\not=l$, $E[\eta
_k^2]=1$, and
\[
\sum_{k=-\infty}^{\infty} a_{nkk}=\sum_{i,j=1}^n K \biggl(\frac{i-j}{nh} \biggr)
\sum_{k=-\infty}^{\infty}\psi_{i-k}\psi_{j-k} =\sum_{i,j=1}^n
\gamma(i-j) b_n(i,j).
\]

By Lemma~\ref{lema.5}, to prove Theorem~\ref{them2.4}, it suffices to show that as $n\to
\infty$,
%
\begin{eqnarray}
\label{ji0}\frac{ (\widetilde Q_{1n}-E\widetilde Q_{1n} )}{\tau_n} &\rightarrow_P& 0, \\
\label{ji1}2 \Vert A\Vert^2 &\equiv&2 \sum_{k,l=-\infty}^{\infty}a_{nkl}^2 \sim
A_{0}^2 \tau_n^2,\\
\label{ji1a}V^2&=&\sum_{k=-\infty}^{\infty}a_{nkk}^2 =  \mathrm{o}( \tau_n^2),\\
\label{ji2}\operatorname{Tr}(A^4)&=&\mathrm{o}(\tau_n^4),
\end{eqnarray}
where $\tau_n=n^{2-\al} h^{3/2-\al}$, and $A_0^2=8 \eta^2A_{\alpha
}$ with
$A_{\alpha}$
as defined in (\ref{ad1}). Indeed, by virtue of (\ref{ji1})--(\ref
{ji2}), it follows from
Lemma~\ref{lema.5} that
\[
(\widetilde{Q}_{n}-E[\widetilde{Q}_{n}])/\tau_n \rightarrow_D A_0 N(0,1).
\nonumber
\]
This, together with (\ref{ji0}), yields Theorem~\ref{them2.4}.

In what follows, we give the proofs of (\ref{ji0})--(\ref{ji2}).
We start with (\ref{ji0}). Recall that $1/2<\al<1$ and $\gamma
(k)\sim\eta|k|^{-\al}$.
By virtue of (\ref{fin15}),
it can be readily seen that
\begin{eqnarray*}
E [ (\widetilde Q_{1n}-E\widetilde Q_{1n} )^2 ]
&=& K^2(0) E \Biggl[\sum_{k=1}^n (e_k^2-Ee_k^2) \Biggr]^2 \no\\
&\le&
C \sum_{j,k=1}^n \gamma^2(j-k)
\le C \sum_{k=1}^n (n-k) \gamma^2(k) \le C n.
\end{eqnarray*}
This, together with the Markov inequality and $nh\to\infty$, yields
(\ref{ji0}).

We next prove (\ref{ji1}). We have
\begin{eqnarray*}\label{ji18}
\Vert A\Vert^2 &=&\sum_{k,l=-\infty}^{\infty} \Biggl(
\sum_{i,j=1}^n\psi_{i-k}\psi_{j-l} K \biggl(\frac{i-j}{nh} \biggr) \Biggr)^2\\
&=&\sum_{i,j,
s,t=1}^n\sum_{k,l=-\infty}^{\infty} \psi_{i-k}\psi_{j-l}\psi
_{s-k}\psi_{t-l}
K\biggl (\frac{i-j}{nh} \biggr) K \biggl(\frac{s-t}{nh} \biggr)\\
&=&\sum_{i,j, s,t=1}^n K \biggl(\frac{i-j}{nh} \biggr) K \biggl(\frac{s-t}{nh} \biggr)
\gamma(i-s) \gamma(t-j).
\end{eqnarray*}
Write
$
f_n(x,y; z,w) =K (\f{x-z}{nh} ) K (\f{y-w}{nh} ) +
K (\f{y-z}{nh} ) K (\f{x-w}{nh} ).
$
Clearly, $f_n(\cdot, \cdot;\cdot,\cdot)$ has the following symmetry
in its
indexes:
\[
f_n(x,y; z,w)=f_n(y,x; z,w)=f_n(x,y; w,z)=f_n(y,x; w,z).
\]
Also $f_n(x,y; z,w)=f_n(z,w; x,y)$. Noting that for any
function $g(x, y)$ and symmetric function $b(x)$,
\[
\sum_{i,j=1}^n
b(i-j)g(i,j)=b(0)\sum_{i=1}^ng(i,i)+\sum_{i=1}^{n-1}\sum_{j=1}^{n-i}
b(i) [g(j, j+i)+g(j+i, j) ],
\]
some algebra shows that (noting that $\gamma(k)=\gamma(-k)$)
%
\begin{eqnarray}
2 \Vert A\Vert^2
&=& \sum_{i_1, i_2,j_1,j_{2}=1}^n
\gamma(i_1-i_2) \gamma(j_1-j_2) f_n(i_1,i_2;j_i,j_2)
\nonumber\\
&=&\gamma^2(0) \sum_{i, j=1}^nf_n(i,i; j,j)
+ 4 \gamma(0) \sum_{i=1}^n\sum_{j=1}^{n-1}\sum_{k=1}^{n-j} \gamma(i)
f_n(i,i; k,j+k)
\label{ji3}\\
&&{} +  4 \sum_{i_1=1}^{n-1} \sum_{i_2=1}^{n-i_1} \sum_{j_1=1}^{n-1}
\sum_{j_2=1}^{n-j_1} \gamma(i_1)
\gamma(j_1) f_n(i_2,i_1+i_2; j_2,j_1+j_2)\no\\
&\equiv& \Delta_{1n} + 4 \Delta_{2n} +4 \Delta_{3n}.
\nonumber
\end{eqnarray}

Recalling that $K(x)$ is a probability density function and
$\gamma(x)\sim\eta x^{-\alpha}$ for $0<\al<1$ and $x>0$, we have that
%
\begin{equation}\label{kai}
\Delta_{3n} \sim
\int_{1}^{n-1} \int_{1}^{n-1} \int_{1}^{n-x}
\int_{1}^{n-y} \gamma(x) \gamma(y) f_n(z,x+z; w,y+w)
\,\mathrm{d}x\,\mathrm{d}y\,\mathrm{d}z\,\mathrm{d}w.
\end{equation}
Write $g_n(x,y)=\int_{1}^{n-x} \int_{1}^{n-y} f_n(z,x+z; w,y+w) \,\mathrm{d}z\,\mathrm{d}w$.
By Fubini's
theorem,
\begin{eqnarray*}
g_n(x,y) &=& \int_{1}^{n-x}
\int_{1-z}^{n-y-z} f_n(z,x+z; w+z,y+w+z) \,\mathrm{d}w \,\mathrm{d}z \\
&=&\int_{-(n-x-1)}^{n-y-1}
\int_{\max\{1, 1-w\}}^{\min\{n-x, n-y-w\}} f_n(0,x; w,y+w) \,\mathrm{d}z \,\mathrm{d}w
\\
&=& \int_{0}^{n-y-1} (n-1-\max\{x, y+w\} ) f_n(0,x; w,y+w) \,\mathrm{d}w
\\
&&{} + \int_{0}^{n-x-1} (n-1-\max\{x+w,
y\} ) f_n(0,x; -w,y-w) \,\mathrm{d}w \no\\
&=& \int_{y}^{n-1} (n-1-\max\{x, w\} ) f_n(0,x; w-y,w) \,\mathrm{d}w \\
&&{} + \int_{x}^{n-1} (n-1-\max\{w, y\} ) f_n(0,x; x-w,x+y-w)
\,\mathrm{d}w.
\end{eqnarray*}

Substituting this into (\ref{kai}), simple calculations show that
if $h\to0$ and $nh\to\infty$, then
%
\begin{eqnarray}\label{wg1}
\Delta_{3n} &\sim&\int_{1}^{n-1} \int_{1}^{n-1} \gamma(x)
\gamma(y) g_n(x,y) \,\mathrm{d}x\,\mathrm{d}y \no\\
&\sim&2 \int_{1}^{n-1} \int_{1}^{n-1}\int_{x}^{n-1}
\gamma(x) \gamma(y) (n-1-\max\{w, y\} ) \no\\
&&\hphantom{\int_{1}^{n-1} \int_{1}^{n-1}\int_{x}^{n-1}}{}
\times\biggl[K \biggl(\f{w}{nh}\biggr )K \biggl(\f{x+y-w}{nh} \biggr){}
+ K \biggl(\f{w-x}{nh} \biggr)K \biggl(\f
{w-y}{nh} \biggr) \biggr] \,\mathrm{d}x\,\mathrm{d}y\,\mathrm{d}w\no\\
&\sim&2 n (nh)^{3-2\al} \eta^2 \int_{0}^{1/h} \int_{0}^{1/h}
\int_{x}^{1/h} x^{-\al} y^{-\al} (1-h\max\{w, y\} )\no\\
&&\hphantom{2 n (nh)^{3-2\al} \eta^2 \int_{0}^{1/h} \int_{0}^{1/h}
\int_{x}^{1/h}}{}\times[K(w)K(x+y-w)\nonumber\\
&&\hphantom{\hphantom{2 n (nh)^{3-2\al} \eta^2 \int_{0}^{1/h} \int_{0}^{1/h}
\int_{x}^{1/h}}{}\times[}{} + K(w-x)K(w-y) ] \,\mathrm{d}x\,\mathrm{d}y\,\mathrm{d}w
\\
&\sim&
2 n (nh)^{3-2\al} \eta^2 \int_{0}^{\infty} \int_{0}^{\infty}
\int_{x}^{\infty} x^{-\al} y^{-\al} \no\\
&&\hphantom{2 n (nh)^{3-2\al} \eta^2 \int_{0}^{\infty} \int_{0}^{\infty}
\int_{x}^{\infty}}{} \times[K(w)K(x+y-w) + K(w-x)K(w-y) ] \,\mathrm{d}x\,\mathrm{d}y\,\mathrm{d}w \no\\
&= &2 n (nh)^{3-2\al} \eta^2 A_{\al}= n (nh)^{3-2\al} A_{0}^2/4,
\nonumber
\end{eqnarray}
based on the fact that $K(x)$ is symmetric, $A_{\alpha}^2<\infty$,
and (\ref{last1}).

By a similar argument, if
$h\to0$ and $nh\to\infty$, then
%
\begin{equation}\Delta_{1n} + 4 \Delta_{2n}
= \mathrm{O}(n^{3-\al}h^2)= \mathrm{o}(\Delta_{3n}). \label{wg2}
\end{equation}
By virtue
of (\ref{ji3}), (\ref{wg1}) and (\ref{wg2}), we obtain the proof
of (\ref{ji1}).

Third, we prove (\ref{ji1a}). Let
\[
h(i,j,s,t)=\sum_{k=-\infty}^{\infty} \psi_{i-k}\psi_{j-k}\psi
_{s-k}\psi_{t-k}
=\sum_{k=-\infty}^{\infty} \psi_{k}\psi_{j-i+k}\psi_{s-i+k}\psi_{t-i+k}.
\]
By $\psi_j\ge0$ and $K(x)\ge0$, it can be readily seen that for any
$j\ge0$,
$s$, and $t$,
\[
\sum_{i=1}^nh(i,j+i,s,t)\le
\sum_{k=-\infty}^{\infty} \psi_{k}\psi_{j+k} \sum_{i=1}^n\psi
_{s-i+k}\psi_{t-i+k}\le
\gamma(j) \gamma(t-s).
\nonumber
\]

Therefore, as in (\ref{ji3})--(\ref{wg1}), it follows from (\ref{last0})
that
%
\begin{eqnarray}\label{jiti}
V^2 &=&\sum_{k=-\infty}^{\infty} \Biggl(
\sum_{i,j=1}^n\psi_{i-k}\psi_{j-k}K \biggl(\frac{i-j}{nh} \biggr) \Biggr)^2
\no\\
&=&\sum_{i,j, s,t=1}^n K \biggl(\frac{i-j}{nh} \biggr) K \biggl(\frac{s-t}{nh}\biggr ) h(i,j,s,t)
\no\\
&\le& K(0)\sum_{i,s,t=1}^n K \biggl(\frac{s-t}{nh} \biggr) h(i,i,s,t) \no\\
& &{}+  2\sum_{j=1}^{n} K \biggl(\frac{j}{nh} \biggr) \sum_{i,s,t=1}^n
K \biggl(\frac{s-t}{nh} \biggr) h(i,j+i,s,t) \\
&\le&
\Biggl[K(0)\gamma(0)+2\sum_{j=1}^{n} K \biggl(\frac{j}{nh} \biggr) \gamma(j)\Biggr ]
\sum_{s,t=1}^n K \biggl(\frac{s-t}{nh} \biggr) \gamma(s-t)\no\\
&\le&C n \biggl(K(0)\gamma(0)+2\int_1^n
x^{-\alpha}K \biggl(\frac{x}{nh} \biggr) \,\mathrm{d}x \biggr)
\int_1^n x^{-\alpha}K \biggl(\frac{x}{nh} \biggr) \,\mathrm{d}x \no\\
&\le& C n^{3-2\alpha} h^{2-2\alpha} =\mathrm{o}(\tau_n^2)
\nonumber
\end{eqnarray}
because $nh\to\infty$. This proves (\ref{ji1a}).

Finally, we prove (\ref{ji2}). Tedious calculations show that
\begin{eqnarray*}
\operatorname{Tr}(A^4)&=&
\sum_{i,j,l,m=-\infty}^{\infty}a_{nij} a_{njl} a_{nlm} a_{nmi}\\
&=& \sum_{i,j,l,m=-\infty}^{\infty} \sum_{j_1,
j_2,\ldots,j_7,j_8=1}^n
\psi_{j_1-i}\psi_{j_2-j}\psi_{j_3-j}\psi_{j_4-l}
\psi_{j_5-l}\psi_{j_6-m}\psi_{j_7-m}\psi_{j_8-i} \\
&&\hphantom{\sum_{i,j,l,m=-\infty}^{\infty} \sum_{j_1,
j_2,\ldots,j_7,j_8=1}^n}{} \times K \biggl(\f{j_1-j_2}{nh} \biggr)
K\biggl (\f{j_{3}-j_4}{nh}\biggr ) K \biggl(\f{j_5-j_6}{nh}\biggr )
K \biggl(\f{j_{7}-j_8}{nh} \biggr)\\
&=&
\sum_{j_1, j_2,\ldots,j_7,j_8=1}^n
K \biggl(\f{j_1-j_2}{nh} \biggr)\gamma(j_2-j_3)\cdots
K \biggl(\f{j_7-j_8}{nh} \biggr)\gamma(j_{8}-j_1).
\end{eqnarray*}
Recall that $K(x)=O[(1+|x|^{1-\beta})^{-1}]$. Similar to the
proof of (\ref{ji1}), it follows from Lemma~\ref{lema.2} that
\begin{eqnarray*}\label{wg3}
\operatorname{Tr}(A^4) &\sim&
\int_{1}^n \int_1^n\cdots\int_{1}^n
K \biggl(\frac{x_1-x_2}{nh} \biggr)\gamma(x_2-x_3)\cdots\\
&&\hphantom{\int_{1}^n \int_1^n\cdots\int_{1}^n}{}\times
K \biggl(\frac{x_{7}-x_8}{nh} \biggr) \gamma(x_8-x_1)
\,\mathrm{d}x_1 \,\mathrm{d} x_2\cdots \,\mathrm{d}x_7\,\mathrm{d}x_{8}\\
&\sim& \eta^4 (nh)^{4(2-\al)}
\int_{0}^{1/h} \int_{0}^{1/h}\cdots\int_{0}^{1/h}
K(x_{1}-x_{2}) |x_{2}-x_{3}|^{-\al}\\
&&\hphantom{\eta^4 (nh)^{4(2-\al)}
\int_{0}^{1/h} \int_{0}^{1/h}\cdots\int_{0}^{1/h}}{}\times K ({x_{7}-x_8} ) |{x_{8}-x_{1}}|^{-\al} \,\mathrm{d}x_1\,\mathrm{d}
x_2\cdots \,\mathrm{d}x_7 \,\mathrm{d}x_{8} \\
&=& \mathrm{O}(1) (nh)^{4(2-\al)}
\int_{0}^{1/h} \int_{0}^{1/h}\cdots\int_{0}^{1/h}
|x_{1}-x_{2}|^{\beta-1} |x_2-x_3|^{-\al} \cdots\\
&&\hphantom{\mathrm{O}(1) (nh)^{4(2-\al)}
\int_{0}^{1/h} \int_{0}^{1/h}\cdots\int_{0}^{1/h}}{}
\times| {x_{7}-x_8}|^{\beta-1} |{x_{8}-x_{1}}|^{-\al}
\,\mathrm{d}x_1\,\mathrm{d} x_2\cdots \,\mathrm{d}x_7
\,\mathrm{d}x_{8}\\
&=&\mathrm{o}(1) (nh)^{4(2-\al)} (1/h)^{2} = \mathrm{o}(1) \tau_n^4.
\end{eqnarray*}
This yields (\ref{ji2}) and thus completes the proof of Theorem~\ref{them2.4}.
\end{pf*}
%
\subsection{\texorpdfstring{Proofs of Theorems \protect\ref{them2.1} and \protect\ref{them2.2}}
{Proofs of Theorems 2.1 and 2.2}}\label{appa.3}
\begin{pf*}{Proof of Theorem~\ref{them2.1}} Note that for any $\alpha$ and
$\widehat{\alpha}$,
%
\begin{eqnarray}\label{l9}
\bigl|n^{(\al-\widetilde\al)}-1 \bigr|\le|\al-\widetilde\al| \log n\cdot
\exp\{|\al-\widetilde\al| \log n\}.
\end{eqnarray}

By (\ref{eq:2.8}), Theorem~\ref{them2.3}, and Assumption~\ref{ass2.2},
the proof of Theorem~\ref{them2.1} follows if we prove
%
\begin{eqnarray}
\label{eq:jiti1}\frac1n \sum_{t=1}^n\widehat f(X_t) &=&
\frac1{n^2h}\sum_{i,j=1}^n K \biggl(\f{X_i-X_j}{h} \biggr) \rightarrow_P
\int_{-\infty}^{\infty} f^2(x)\,\mathrm{d}x,
\\
\label{eq:jiti2}\frac1n \sum_{t=1}^n \widehat e_t^2 &=& \frac1n \sum_{t=1}^n
\bigl(e_t+m_{\theta_0}(X_t)-m_{\widetilde\theta}(X_t)\bigr)^2\rightarrow_P
\gamma(0),
\end{eqnarray}
and under the relevant conditions of Theorem~\ref{them2.1},
%
\begin{equation}\label{eq:jiti3}
2 R_{1n}(h)+R_{2n}(h)  =  \mathrm{o}_P(\sigma_{in}(h)),\qquad i=1, 2,
\end{equation}
where $R_{1n}(h)$ and $R_{2n}(h)$ are defined as in (\ref{eq:2.8}).

Recall that $E [K (\f{X_1-X_2}{h} ) ]\sim h \int_{-\infty}^{\infty}
f^2(x)\,\mathrm{d}x$ by (\ref{try1}).
The proof of (\ref{eq:jiti1}) follows from standard methodologym and
thus the details are omitted.
By the stationary ergodic theorem, $\frac1n\sum_{t=1}^ne_t^2\to
\gamma(0)$,
a.s.
This implies that (\ref{eq:jiti2}) will follow if we have
%
\begin{equation}\label{f2}
\frac1n \sum_{t=1}^n \bigl[
2\ep_t\bigl(m_{\theta_0}(X_t)-m_{\tilde\theta}(X_t)\bigr)+
\bigl(m_{\theta_0}(X_t)-m_{\tilde\theta}(X_t)\bigr)^2 \bigr]\rightarrow_P 0.
\end{equation}
Because the proof of (\ref{f2}) is similar to that of (\ref
{eq:jiti3}), we only prove (\ref{eq:jiti3}) what follows.

For $\forall\ep>0$, write $\Omega_n=\{\widetilde{\theta}\dvt
\Vert\widetilde{\theta}-\theta_0\Vert\le\ep n^{-\alpha/2}\}$, and let
\begin{eqnarray*}
J_1(s,t) &=& K \biggl(\frac{X_s - X_t}{h} \biggr) \frac{\partial m_{\theta
}(X_s)}{\partial
\theta}\bigg|_{\theta=\theta_0}, \\
J_2(s,t) &=& K \biggl(\frac{X_s - X_t}{h} \biggr) \biggl\{m_{\theta_0}(X_s) -
m_{\widetilde{\theta}}(X_s) + (\theta_0 -
\widetilde{\theta} )^{\tau} \frac{\partial m_{\theta
}(X_s)}{\partial
\theta}\bigg|_{\theta=\theta_0} \biggr\}
\end{eqnarray*}
and
\[
J_3(s,t) = K \biggl(\frac{X_s - X_t}{h} \biggr) \bigl(m_{\theta_0}(X_s) -
m_{\widetilde{\theta}}(X_s) \bigr) \bigl(m_{\theta_0}(X_t) -
m_{\widetilde{\theta}}(X_t) \bigr).
\label{eq:gao1}
\]
Using this notation, we have $R_{2n}(h) = \sum_{t=1}^n \sum_{s=1,\not=t}^n
J_3(s,t)$ and
\[
R_{1n} (h)= (\theta_0 -\widetilde{\theta} ) \sum_{t=1}^n e_t\sum
_{s=1,\not=t}^n
J_1(s,t)+ \sum_{t=1}^n e_t\sum_{s=1,\not=t}^n J_2(s,t).
\]

Recalling Assumptions~\ref{ass2.2}, \ref{ass2.4}(i) and~\ref{ass2.5}, it is readily apparent that
\begin{eqnarray*}\label{eq:gao2}
E|J_1(s, t)| &\le&
E \biggl[K \biggl(\frac{X_s - X_t}{h} \biggr)
\biggl\Vert\frac{\partial m_{\theta}(X_s)}{\partial
\theta}\bigg|_{\theta=\theta_0}
\biggr\Vert \biggr]\\
&\le& C h
E
\biggl\Vert\frac{\partial m_{\theta}(X_s)}{\partial
\theta}\bigg|_{\theta=\theta_0}
\biggr\Vert
\leq C_1 h,\\
E [|J_1(s, t)J_1(s_1, t)| ]
&\le&
E \biggl[K \biggl(\frac{X_s - X_t}{h} \biggr)K \biggl(\frac{X_{s_1} -
X_t}{h} \biggr)\\
&& \hphantom{E \biggl[}{}\times\biggl\Vert\frac{\partial m_{\theta}(X_{s})}{\partial
\theta}\bigg|_{\theta=\theta_0}
\biggr\Vert \biggl\Vert\frac{\partial m_{\theta}(X_{s_1})}{\partial
\theta}\bigg|_{\theta=\theta_0}
\biggr\Vert \biggr]\\
&\le& C h^2
E \biggl[ \biggl\Vert\frac{\partial m_{\theta}(X_1)}{\partial
\theta}\bigg|_{\theta=\theta_0}
\biggr\Vert^2 \biggr]
\leq C_1 h^2
\end{eqnarray*}
for all different values of $s, s_1$, and $t$.
These findings imply that for any $1\le t\le n$,
\[\label{eq:gao3a}
E \Biggl[\sum_{s=1,\not=t}^n
J_1(s,t) \Biggr]^2\le
C (nh+n^2h^2)
\le2C (nh)^2,
\]
because $nh\to\infty$. Thus, by the independence of $e_t$ and $X_s$,
\begingroup
\abovedisplayskip=7pt
\belowdisplayskip=7pt
\begin{equation}\label{eq:gao4}
 E \Biggl[\sum_{t=1}^n e_t\sum_{s=1,\not=t}^n
J_1(s,t)\Biggr ]^2 \le C (nh)^2 \sum_{t_1, t_2=1}^n
E [e_{t_1}e_{t_2} ]
\le C_1 n^{4-\alpha}h^2,
\end{equation}
where we have used the fact that
\[
\sum_{t_1, t_2=1}^n E [e_{t_1} e_{t_2} ]= E \Biggl(\sum_{t=1}^n
e_t\Biggr )^2\le Cn^{2-\alpha},
\]
as seen from (\ref{ad10}).
On the other hand, it follows from Taylor's expansion of
$m_{\theta}(x)$ (with respect to $\theta$) that
under $H_0$, for all $s\neq s_1$, $s\neq t$ and $s_1\neq t$ and for $n$
large enough such
that $\Omega_n\subseteq\Theta_0$,
%
\begin{eqnarray}\label{lp19}
&& E [|e_t||J_2(s, t)|I{(\widetilde{\theta}\in\Omega_n)} ] \no\\
&&\qquad \le C \ep n^{-\alpha} E|e_t|
E \biggl[K \biggl(\frac{X_s - X_t}{h} \biggr)
\biggl\Vert\frac{\partial^2 m_{\theta}(X_s)}{\partial
\theta^2}\bigg|_{\theta=\theta_0}
\biggr\Vert \biggr]
\\
&&\qquad \le C_1 \ep n^{-\alpha}h E \biggl[
\biggl\Vert\frac{\partial^2 m_{\theta}(X_s)}{\partial
\theta^2}\bigg|_{\theta=\theta_0}
\biggr\Vert^2 \biggr]
\leq C_2 \ep n^{-\alpha}h. \nonumber
\end{eqnarray}

It follows from (\ref{eq:gao4}) and (\ref{lp19}) that
\begin{eqnarray}\label{lp20}
E [|R_{1n}(h)|I{(\widetilde{\theta}\in\Omega_n)} ] &\le& \ep
n^{-\al/2} E \Biggl|\sum_{t=1}^n e_t\sum_{s=1,\not=t}^n
J_1(s,t) \Biggr|\no\\
&& {}+ \sum_{t=1}^n \sum_{s=1,\not=t}^n
E [|e_t| |J_2(s,t)I{(\widetilde{\theta}\in\Omega_n)} ] \\
&\le& C \ep n^{2-\al} h.\no
\end{eqnarray}
This, together with Markov's inequality, yields the following results
for $\forall\ep>0$ and $n$ sufficiently large:
(i) If $nh\to\infty$ and $n^{2(1-\alpha)}h\to0$, then
%
\begin{eqnarray}\label{eq:gao5}
P \bigl( |R_{1n}(h)|\ge\ep^{1/2}\si_{1n} \bigr) &\le&
P (\Vert\widetilde{\theta}-\theta_0\Vert> \ep n^{-\alpha/2} )\no\\
&&{} + C\ep^{-1/2}(n^2h)^{-1/2}E [ | R_{1n}(h) |I{(\widetilde{\theta
}\in
\Omega_n)} ]
\\
&\le& P (\Vert\widetilde\theta-\theta_0\Vert> \ep n^{-\alpha/2} )+
Cn^{1-\alpha}h^{1/2}\ep^{1/2} \le C_1 \ep.\no
\end{eqnarray}

(ii) If $h\to0$ and $n^{2(1-\alpha)}h\to\infty$, then
%
\begin{eqnarray}
\label{eq:gao6}
P \bigl( |R_{1n}(h)|\ge\ep^{1/2}\si_{2n} \bigr) &\le&
P (\Vert\widetilde{\theta}-\theta_0\Vert> \ep n^{-\alpha/2} )
\nonumber\\
& &{} +
C\ep^{-1/2}(n^{4-2\alpha}h^2)^{-1/2}E [ |R_{1n}(h) |I{(\widetilde
{\theta} \in
\Omega_n)} ]\\
&\le& P (\Vert\widetilde{\theta}-\theta_0\Vert> \ep n^{-\alpha/2} )+
C \ep^{1/2} \le C_1 \ep^{1/2}.\nonumber
\end{eqnarray}
\endgroup
The results (\ref{eq:gao5}) and (\ref{eq:gao6}) yield that
$R_{1n}(h)=\mathrm{o}(\si_{jn})$,
$1\leq j\leq2$, under the corresponding conditions of
Theorems~\ref{them2.1}.
Similarly, by noting
\begin{eqnarray*}\label{eq:gao7}
E [|J_3(s, t)|I{(\widetilde{\theta} \in\Omega_n)} ] &\le&
C \ep n^{-\alpha}
E \biggl[K \biggl(\frac{X_s - X_t}{h}\biggr )
\biggl\Vert\frac{\partial m_{\theta}(X_s)}{\partial
\theta}\bigg|_{\theta=\theta_0}
\biggr\Vert
\biggl\Vert\frac{\partial m_{\theta}(X_t)}
{\partial\theta}\bigg|_{\theta=\theta_0}
\biggr\Vert \biggr]
\\
&\le& C_1 \ep n^{-\alpha}h E\biggl [ \biggl\Vert\frac{\partial m_{\theta
}(X_1)}{\partial
\theta}\bigg|_{\theta=\theta_0}
\biggr\Vert^2 \biggr] \le C_2 \ep n^{-\alpha} h,
\end{eqnarray*}
we obtain that for $1\leq j\leq2$ and $n$ sufficiently large,
%
\begin{eqnarray}\label{eq:gao8}
P\bigl( |R_{2n}(h)|\ge\ep^{1/2}\si_{jn}\bigr) &\le&
P (\Vert\widetilde{\theta}-\theta_0\Vert> \ep n^{-\alpha/2} )
\nonumber
\\[-8pt]
\\[-8pt]
&&{} + C\ep^{1/2}(\si_{jn})^{-1}n^2 E [|J_3(1,
2)|I{(\widetilde{\theta}\in\Omega_n)} ]
\leq C \ep^{1/2},\nonumber
\end{eqnarray}
implying that $R_{2n}(h)=\mathrm{o}(\si_{jn})$ holds for $1\leq j\leq2$.
Combining these findings, we obtain (\ref{eq:jiti3}), and thus also
complete the proof of Theorem~\ref{them2.1}.
\end{pf*}
\begin{pf*}{Proof of Theorem~\ref{them2.2}}
 As in (\ref{eq:2.8}), under $H_0$, we can
write
%
\begin{equation}\label{eq:2.8a}
 \sum_{t=1}^n \sum_{s=1,\neq t}^n b_n(s,t)
\bigl(\widehat{e}_s \widehat{e}_t - \widehat{\gamma}(s-t) \bigr)
= M_n^*(h) + 2 R_{1n}^*(h)+R_{2n}^*(h)+R_{3n}(h),
\end{equation}
where $M_n^*(h) = \sum_{t=1}^n \sum_{s=1,\neq t}^n b_n(s, t) [e_s
e_t-\gamma(s-t) ] $,
\begin{eqnarray*}
R_{1n}^*(h)
&=&
\sum_{t=1}^n \sum_{s=1,\neq t}^n b_n(s, t) e_s
[m_{\theta_0} (t/n )
- m_{\widetilde{\theta}} ( t/n ) ] \\
&=&(\widetilde\theta-\theta_0) \sum_{t=1}^n \sum_{s=1,\neq t}^n
b_n(s, t) e_s
\frac{\partial m_{\theta}(t/n)}{\partial
\theta}\biggm|_{\theta=\theta_0} \\
&&{} +  \sum_{t=1}^n \sum_{s=1,\neq t}^n b_n(s, t) e_s
\biggl[m_{\theta_0} (t/n ) - m_{\widetilde{\theta}} ( t/n ) +(\widetilde
\theta-\theta_0)^{\tau} \frac{\partial m_{\theta}(t/n)}{\partial
\theta}\biggm|_{\theta=\theta_0} \biggr], \\
R_{2n}^*(h) &=&
\sum_{t=1}^n \sum_{s=1,\neq t}^n b_n(s, t) [m_{\theta_0} ( s/n )
- m_{\widetilde{\theta}} ( s/n ) ] [m_{\theta_0} ( t/n )
- m_{\widetilde{\theta}} ( t/n ) ],
\end{eqnarray*}
and by the symmetries of $K(x)$, $\gamma(k)$, and $\hat\gamma(k)$,
%
\begin{equation}\label{ad90}
\hspace*{-8.5pt}R_{3n}(h) = \sum_{t=1}^n \sum_{s=1,\neq t}^n\!\! b_n(s,t)
[ \widehat{\gamma}(s-t)-\gamma(s-t) ]
= 2 \sum_{s=1 }^{n-1}(n-s)K \biggl(\f s{nh}\biggr )
[ \widehat{\gamma}(s)-\gamma(s) ].
\end{equation}
By (\ref{eq:2.8a}) and Theorem
\ref{them2.4}, Theorem~\ref{them2.2} will follow if we prove that
%
\begin{eqnarray}
\label{k6}2 R_{1n}^*(h)+R_{2n}^*(h)&=&\mathrm{o}_P [\sigma_{3n}(h) ], \\[-2pt]
\label{k7}R_{3n}(h)&=&\mathrm{o}_P [\sigma_{3n}(h) ], \\[-2pt]
\label{k8}\f{\widehat\sigma_{3n}(h)}{\sigma_{3n}(h)}&\rightarrow_P&1.
\end{eqnarray}

The proofs of (\ref{k6})--(\ref{k8}) are quite technical and thus are
omitted here, but can be derived from the rest of the proof of Theorem
2.2 in Gao and Wang \cite{gawa10}.
\end{pf*}
%
\section{\texorpdfstring{Proofs of Theorems \protect\ref{them4.1} and \protect\ref{them4.2}}
{Proofs of Theorems 4.1 and 4.2}}\label{appb}
\begin{pf*}{Proof of Theorem~\ref{them4.1}}
We first prove (\ref{eq:3.1.1}). In view of Theorem~\ref{them2.2}, it suffices
to show
that
%
\begin{equation}\label{eq:3.1a}
\sup_{x \in R} \bigl|P^{\ast}\bigl(\widehat{T}_n^{\ast}(h)\leq x\bigr) - \Phi(x)
\bigr| =
\mathrm{o}_P(1).
\end{equation}

As in (\ref{eq:2.8a}), we may rewrite $\widehat{T}_n^{\ast}(h)$ as
%
\begin{equation}
\label{eq:2.8ab}
\widehat{T}_n^{\ast}(h)= \frac1{\widehat\sigma_{3n}(h)} [
M_n^{**}(h) + 2 R_{1n}^{**}(h)+R_{2n}^{**}(h)+R_{3n}^{**}(h) ],
\end{equation}
where $M_n^{**}(h) = \sum_{t=1}^n \sum_{s=1,\neq t}^n b_n(s,
t) [e_s^*
e_t^*-\gamma_{\widetilde\lam}(s-t) ]$,
\begin{eqnarray*}
R_{1n}^{**}(h)
&=& \sum_{t=1}^n \sum_{s=1,\neq t}^n b_n(s, t) e_s^*
[m_{\widetilde\theta} (t/n )
- m_{\widetilde{\theta}^*} ( t/n ) ],
\\[-2pt]
&=&(\widetilde\theta-\widetilde\theta^*) \sum_{t=1}^n \sum
_{s=1,\neq t}^n b_n(s, t) e_s^*
\frac{\partial m_{\theta}(t/n)}{\partial
\theta}\bigg|_{\theta=\widetilde\theta}
\\[-2pt]
& &{}+  \sum_{t=1}^n \sum_{s=1,\neq t}^n b_n(s, t) e_s^*
\biggl[m_{\widetilde\theta} (t/n ) - m_{\widetilde{\theta}^*} ( t/n )
+(\widetilde\theta^*-\widetilde\theta)^{\tau} \frac{\partial
m_{\theta}(t/n)}{\partial\theta}\biggm|_{\theta=\widetilde\theta}\biggr ],
\\[-2pt]
R_{2n}^{**}(h) &=&
\sum_{t=1}^n \sum_{s=1,\neq t}^n b_n(s, t)
\bigl(m_{\widetilde{\theta}^*}(s/n) - m_{\widetilde{\theta}}(s/n)\bigr )
\bigl(m_{\widetilde{\theta}^*}(t/n)
- m_{\widetilde{\theta}}(t/n)\bigr ),
\\[-2pt]
R_{3n}^{**}(h) &=&\sum_{t=1}^n \sum_{s=1,\neq t}^n b_n(s, t)
[\widehat\gamma^*(s-t)-\gamma_{\widetilde\lam}(s-t) ],
\end{eqnarray*}
where $\widehat{e}_s^{\,*}=Y_s^*-m_{\tilde\theta^*}(X_s)$ and
\[
\widehat\gamma^*(k) = \cases{
\displaystyle\frac1n \sum_{i=1}^{n-|k|}\widehat{e}_i^{\,*}\widehat{e}_{i+|k|}^{\,*}
& \quad for  $|k|\le(nh)^{1/3}$,\cr
\widetilde{\eta} |k|^{-\widetilde\al} &\quad for $(nh)^{1/3}<|k|\le
n-1$.}
\]

The result (\ref{eq:3.1a}) will follow if we prove that
%
\begin{equation}\label{o2}
I_{(\tilde\al\in\Lambda_n)} \sup_{x \in R} \biggl|P \biggl(\frac
{M_n^{**}(h)}{\widehat\sigma_{3n}(h)}\leq
x\Big\vert\mathcal{W}_n \biggr)
- \Phi(x) \biggr| = \mathrm{o}_P(1),
\end{equation}
and, for any $\ep>0$,
%
\begin{equation}\label{o3}
I_{\{(\tilde\al\in\Lambda_n, \tilde\theta\in\Omega_{1n}\}} P
\bigl(|2 R_{1n}^{**}(h)+R_{2n}^{**}(h)+R_{3n}^{**}(h)|\ge\ep\widehat
\sigma_{3n}(h)\vert\mathcal{W}_n \bigr) = \mathrm{o}_P(1),
\end{equation}
where $\mathcal{W}_n=(Y_1, \ldots,
Y_n)$, $\Lambda_n=\{\tilde\al\dvt |\tilde\al-\al|\le C
{w}_n^{-1/2}\}$
and
$
\Omega_{1n} =\{\tilde\theta\dvt \Vert\tilde\theta-\theta_0\Vert\le
C_0n^{-\al/2}\}
$,
with $C_0$ chosen so that $P(\Vert\tilde\theta-\theta_0\Vert\ge C_0n^{-\al
/2})\le\ep$.
Indeed, recalling Assumptions~\ref{ass2.3} and~\ref{ass2.6}(ii),
we have $I(\tilde\al\notin\Lambda_n \mbox{ or } \tilde\theta
\notin\Omega_{1n})=\mathrm{o}_P(1)$. This, together with
(\ref{o2}) and (\ref{o3}), proves the statement.

We next prove (\ref{o2}) and (\ref{o3}). We start with (\ref{o2}).
As constructed in Section~\ref{sec4.1}, $\{e_t^{\ast}\}$ is written as
%
\begin{equation}
\label{111}
e_t^{\ast} =\sum_{j=-\infty}^{\infty} \psi^*_j \eta_{t-j}^{\ast}
\end{equation}
such that $E[e_t^{\ast}e_{t+k}^{\ast}|\mathcal{W}_n]={\gamma
}_{\widetilde{\lambda}}(k)$.
Also note that $\hat\al\in\Lambda_n$ and $1/2<\al<1$ imply that
there exists a $\delta_0>0$
such that $1/2+\delta_0<\tilde\al<1-\delta_0$ whenever $n$ is
sufficiently large.
Based on these facts,
given $\mathcal{W}_n$, $M_n^{**}(h)$ has the same structure as that of
$Q_{1n}-E[Q_{1n}]$
defined in the proof of
Theorem~\ref{them2.4}. It follows from the same
routine as in the proof of Theorem~\ref{them2.4} that
%
\begin{equation}\label{tyu1}
I_{(\tilde\al\in\Lambda_n)} \sup_{x \in R} \bigl|P \bigl( {M_n^{**}(h)}\leq
{\widehat\tau_n\widehat A_0} x\vert\mathcal{W}_n \bigr)
- \Phi(x) \bigr| = \mathrm{o}_P(1),
\end{equation}
where $\widehat\tau_n=n^{2-\tilde\al}h^{3/2-\tilde\al}$ and $\hat
A_0^2=8\tilde\eta^2A_{\tilde\al}$
with
\begin{eqnarray*}
A_{\tilde\alpha}&=&\int_0^{\infty} \int_0^{\infty} \int
_0^{\infty} x^{-\tilde\alpha}
y^{-\tilde\alpha}
[K(z)K(x+y-z) + K(z-x)K(z-y) ] \,\mathrm{d}x \,\mathrm{d}y \,\mathrm{d}z.
\end{eqnarray*}
Recall that Assumption~\ref{ass2.4}(ii) implies that for any $u\in\R$,
%
\begin{equation}\label{last3}
\int_{0}^{\infty} K(w)K(w+u) \,\mathrm{d}w \le C/(1+|u|^{1-\beta}).
\end{equation}
It now
can be readily seen that whenever $1/2+\delta_0<\tilde\al<1-\delta_0$,
\begin{eqnarray*}
A_{\tilde\alpha}&=&\int_{1/n}^{n} \int_{1/n}^{n} \int_{1/n}^{n}
x^{-\tilde\alpha}
y^{-\tilde\alpha}
[K(z)K(x+y-z) +  K(z-x)K(z-y) ] \,\mathrm{d}x \,\mathrm{d}y \,\mathrm{d}z+\mathrm{o}_{P}(1)
\\[-2pt]
& =& A_{\tilde\alpha}^* [1+\mathrm{o}_{P}(1)],
\end{eqnarray*}
that is, ${\widehat\tau_n\widehat A_0}=\widehat\sigma_{3n}(h)
[1+\mathrm{o}_{P}(1)]$, where we have used the fact
that $A_{\tilde\alpha}^*-A_{\al}=\mathrm{o}_P(1)$ and $0<A_{\al}<\infty$.
Substituting this into (\ref{tyu1}), we get (\ref{o2}).

The proof of (\ref{o3}) follows from the same arguments as in the
proofs of (\ref{k6}) and (\ref{k7}). Details are given in the rest of
the proof of the first part of Theorem 4.1(i) of Gao and Wang~\cite{gawa10}.

We next prove the second part of Theorem~\ref{them4.1}(i). In view of Theorem
\ref{them2.2}, it suffices to show
that
%
\begin{equation}\label{k10}
l_r^*-l_r=\mathrm{o}_P(1).
\end{equation}
In fact, recalling the definitions of $l_r^*$ and $l_r$, it can be readily
seen
from (\ref{eq:3.1a}) and Theorem~\ref{them2.2} that
$
\Phi(l_r^*)-\Phi(l_r)=\mathrm{o}_P(1),
$
which implies (\ref{k10}), because $\Phi(x)$ is a bounded continuous function.

Finally, we prove Theorem~\ref{them4.1}(ii). In view of (\ref{k10}),
it suffices to show that under $H_1$,
%
\begin{equation}\label{k11}
P\bigl(\widehat L_{3n}(h)\ge l_r\bigr)=1,
\end{equation}
with $l_r$ satisfying $\Phi(l_r)=1-r+\mathrm{o}(1)$ with $0<r<1$. To prove
(\ref{k11}),
as in (\ref{eq:2.8a}), under $H_1$, we may rewrite $\widehat
L_{3n}(h)$ as
%
\begin{equation}
\widehat L_{3n}(h)= \frac1{\widehat\sigma_{3n}(h)} [
S_n(h) +2Q_{1n}(h)+Q_{2n}(h) +R_{3n}(h) ],
\end{equation}
where
\begin{eqnarray*}
S_n(h) &=& \sum_{t=1}^n \sum_{s=1,\neq t}^n b_n(s, t) [\zeta_s
\zeta_t-\gamma(s-t) ],
\\
Q_{1n}(h)&=& c_n \sum_{t=1}^n \sum_{s=1,\neq t}^n b_n(s, t) \zeta_s
\Delta\biggl(\f tn \biggr),
\\
Q_{2n}(h) &=& c_n^2
\sum_{t=1}^n \sum_{s=1,\neq t}^n b_n(s, t) \Delta\biggl(\f sn \biggr) \Delta\biggl(\f
tn \biggr),
\end{eqnarray*}
in which $\zeta_t=e_t+ [m_{\theta_1} (\f tn ) - m_{\widetilde{\theta
}} (\f tn ) ]$ and $R_{3n}(h) $ is defined as in (\ref{eq:2.8a}).
Simple calculations show that
%
\begin{eqnarray}
Q_{2n}(h)&\sim& c_n^2 \int_1^n\int_1^n K \biggl(\f{x-y}{nh}\biggr )\Delta\biggl(\f
xn \biggr)\Delta\biggl(\f yn \biggr)
\,\mathrm{d}x\,\mathrm{d}y \no
\\[-8pt]
\\[-8pt]
&\sim& A_0 [1+\mathrm{o}(1)] c_n^2 (nh)^2=A_0 d_n \sigma_{3n}(h),\no
\end{eqnarray}
where $d_n=c_n^2n^{\al} h^{\al-1/2}\to\infty$ and $A_0>0$.
Recalling Assumption~\ref{ass4.1}(i), the same arguments as in the proof of
Theorem~\ref{them2.2} yield that
%
\begin{equation}
S_n(h)/\sigma_{3n}(h)\rightarrow_D N(0,1),
\end{equation}
and $Q_{1n}(h)=\mathrm{O}_P (c_nn^{-\al/2}(nh)^2 )=\mathrm{o}_P (Q_{2n}(h) )$. These
findings, together with (\ref{k7}) and (\ref{k8}) (i.e., $R_{3n}=\mathrm{o}_P
(\si_{3n}(h) )$ and $\hat\si_{3n}(h)/\si_{3n}(h)\rightarrow_P 1$), imply that
\[
\widehat{ L}_{3n}(h)-A_0d_n^2\rightarrow_D N(0,1),
\]
with $d_n^2\to\infty$ and $A_0>0$. We now have (\ref{k11}), because
$l_r$ is finite for $0<r<1$.
The proof of Theorem~\ref{them4.1} is completed.
\end{pf*}
\begin{pf*}{Proof of Theorem~\ref{them4.2}}
Observe that under $H_0$,
\begingroup
\abovedisplayskip=6.7pt
\belowdisplayskip=6.7pt
\begin{eqnarray}
\label{eq:jitigao1}
u_t & \equiv& \widehat{e}_t = Y_t - m_{\widetilde{\theta}}(X_t) =
e_t + m_{\theta_0}(X_t) - m_{\widetilde{\theta}}(X_t)
\nonumber
\\[-9pt]
\\[-9pt]
& = & e_t + (\theta_0 - \widetilde{\theta})^{\tau} \frac{\partial
m_{\theta}(X_t)}{\partial\theta}\bigg|_{\theta=\theta_0} + \mathrm{o}_P
(\Vert\theta_0 - \widetilde{\theta}\Vert ),\nonumber
\end{eqnarray}
using Assumption~\ref{ass2.6}.

By straightforward calculations, we then have, for $n$ large enough,
%
\begin{eqnarray}
\label{eq:jitigao2}
I_u(\omega) & = & \frac{1}{2\uppi n} \Biggl|\sum_{s=1}^n u_s \mathrm{e}^{\mathrm{i}s \omega} \Biggr|^2
\nonumber
\\[-2pt]
& = & \frac{1}{2\uppi n} \Biggl(\sum_{s=1}^n u_s^2 + 2 \sum_{s=1}^{n-1} \sum
_{t=1}^{n-s} \cos(\omega s) u_{s+t} u_s \Biggr)
\nonumber
\\[-9pt]
\\[-9pt]
& = & \frac{1}{2\uppi n} \Biggl(\sum_{s=1}^n e_s^2 + 2 \sum_{s=1}^{n-1} \sum
_{t=1}^{n-s} \cos(\omega s) e_{s+t} e_s \Biggr) + \mathrm{o}_P(1)
\nonumber
\\[-2pt]
& = & I_e(\omega) + \mathrm{o}_P(1)\nonumber
\end{eqnarray}
using Assumptions~\ref{ass2.1} and~\ref{ass2.6}.

Thus, for $n$ large enough,
%
\begin{eqnarray}
\label{eq:jitigao3}
\Gamma_u(\lambda) & = & \frac{1}{4\uppi} \int_{-\uppi}^{\uppi} \biggl(\log
(\psi(\omega; \lambda)) + \frac{I_u(\omega)}{\psi(\omega; \lambda
)} \biggr) \,\mathrm{d}\omega
\nonumber
\\[-2pt]
& = & \frac{1}{4\uppi} \int_{-\uppi}^{\uppi} \biggl(\log(\psi(\omega;
\lambda)) + \frac{I_e(\omega)}{\psi(\omega; \lambda)} \biggr)\,\mathrm{d}\omega + \mathrm{o}_P(1)
\\[-2pt]
& = & \Gamma_e(\lambda) + \mathrm{o}_P(1).\nonumber
\end{eqnarray}

Therefore, by Assumption 4.3 and Theorem 1(ii) of Heyde and Gay \cite{hega93},
we have, as $n\rightarrow\infty$,
%
\begin{equation}
\label{eq:jitigao4}
\sqrt{n} (\widetilde{\lambda} - \lambda) \rightarrow_D N (0,
\Sigma^{-1}(\lambda) ),
\end{equation}
which implies for $n$ large enough,
%
\begin{equation}
\label{eq:jitigao5}
\frac{n^{2/5}}{\log(n)} (\widetilde{\lambda} - \lambda) \sim\frac
{1}{n^{1/10} \log(n)} \sqrt{n} (\widetilde{\lambda} - \lambda)
\rightarrow_P 0.
\end{equation}

This completes the proof of Theorem~\ref{them4.2}.
\end{pf*}
\endgroup
\end{appendix}\eject
\section*{Acknowledgements} The authors thank the Editor, the
Associate Editor, and a referee for their constructive comments that
have led to a significant improvements in the original manuscript.
Earlier versions of this paper were presented at several conferences
and seminars. The authors thank the seminar and conference
participants, in particular Vo Anh, Maxwell King, Nicolai Leonenko,
Zhengyan Lin, Eckhard Platen, and Qiman Shao, for constructive
suggestions and comments. Thanks also go to the Australian Research
Council Discovery Grants Program for its financial support.


\printhistory

\end{document}